\documentclass{elsart}

\usepackage{amssymb}

\usepackage{amsmath}
\usepackage{enumerate}
\usepackage{url}

\DeclareMathOperator{\Bad}{\bf{Bad}}
\newcommand{\cthbn}{\ensuremath{\mathcal{C}(\theta B_n)}}
\newcommand{\cthfn}{\ensuremath{\mathcal{C}(\theta F_n)}}

\theoremstyle{plain}

\begin{document}

\begin{frontmatter}

\title{Diophantine approximation and  badly approximable sets}

\author[simon]{Simon Kristensen\thanksref{label1}}
\ead{Simon.Kristensen@ed.ac.uk}
\author[becky]{Rebecca Thorn}
\ead{r.e.thorn@qmul.ac.uk}
\author[sanju]{Sanju Velani\corauthref{cor1}\thanksref{label2}}
\ead{slv3@york.ac.uk}

\address[simon]{School of Mathematics, The University of Edinburgh,
  James Clerk Maxwell Building, King's Buildings, Mayfield Road,
  Edinburgh EH9 3JZ, UK}
\address[becky]{School of Mathematical Sciences, Queen Mary,
  University of London, \\ Mile End Road, London E1 4NS, UK}
\address[sanju]{Department of Mathematics, University of York,
  Heslington, \\  York YO10 5DD, UK}

\corauth[cor1]{Corresponding author}

\thanks[label1]{William Gordon Seggie Brown Fellow}
\thanks[label2]{Royal Society University Research Fellow}


\begin{abstract}
Let $(X,d)$ be a metric space and $(\Omega, d)$ a
compact subspace of $X$ which supports a non-atomic finite measure
$m$.  We consider `natural' classes of badly approximable  subsets
of $\Omega$. Loosely speaking, these consist of points in $\Omega$
which   `stay clear' of some given set of points in $X$. The
classical set $\Bad$ of `badly approximable' numbers in the theory
of Diophantine approximation falls within our framework as do the
sets $\Bad(i,j)$ of simultaneously badly approximable numbers.
Under various natural conditions we prove that the  badly
approximable subsets of $\Omega$ have full Hausdorff dimension.
Applications of our general framework include those from number
theory (classical, complex, $p$-adic and formal power series)  and
dynamical systems (iterated function schemes, rational maps and
Kleinian groups).
\end{abstract}

\begin{keyword}
Diophantine approximation, Hausdorff dimension, badly approximable
elements, dynamical systems.

\vspace{1ex}
\MSC 11J83 \sep 11J61 \sep 28A80 \sep 37C45 \sep 37F
\end{keyword}

\dedicated{Dedicated to Lalji and Manchaben on their 70 plus
birthdays}

\end{frontmatter}

\section{Introduction}
 \label{intro}

\subsection{ The setup and the problem}
\label{setup}

Let $(X,d)$ be a metric space and $(\Omega, d)$ a compact subspace of
$X$ which contains the support of a non-atomic finite measure $m$.
Let $\mathcal{R}=\{R_\alpha \in X : \alpha \in J\}$ be a family of
subsets $R_\alpha$ of $X$ indexed by an infinite, countable set $J$.
The sets $R_\alpha$ will be referred to as \emph{resonant sets}. Next,
let $\beta: J \rightarrow \mathbb{R}^+ : \alpha \rightarrow
\beta_\alpha$ be a positive function on $J$. To avoid pathological
situations within our framework, we shall assume that the number of
$\alpha \in J$ with $\beta_\alpha$ bounded above is finite -- thus
$\beta_\alpha $ tends to infinity as $\alpha$ runs through $J$.  Given
a real, positive function $\rho : \mathbb{R}^+ \rightarrow
\mathbb{R}^+ : r \to \rho(r) $ such that $\rho(r) \rightarrow 0$ as $r
\rightarrow \infty$ and that $\rho$ is decreasing for $r$ large
enough, consider the set
\begin{multline*}
  \Bad^*(\mathcal{R}, \beta, \rho):=\{x \in \Omega :\exists \; c(x)>0
  \text{ such that } \; d(x,R_\alpha) \geq c(x)\rho(\beta_\alpha) \\
  \text{ for all } \alpha \in J\} \ ,
\end{multline*}
where $ d(x,R_\alpha) := \inf_{a \in R_\alpha} d(x,a) $.  Loosely
speaking, in the case that the resonant sets are points,
$\Bad^*(\mathcal{R}, \beta, \rho)$ consists of points in $\Omega$
which `stay clear' of `$\rho$-balls' centred at resonant points.
Notice that since the number of $\alpha \in J$ with $\beta_\alpha$
bounded above is finite and $\rho$ is eventually decreasing, the
number of $\alpha \in J$ with $\rho(\beta_\alpha) \geq \varepsilon > 0
$ is finite. In view of this, without loss of generality we shall
assume that the $\sup_{\alpha\in J} \rho(\beta_\alpha)$ is finite.
Otherwise, if $\rho(\beta_\alpha)$ can get arbitrarily large, then
trivially $\Bad^*(\mathcal{R}, \beta, \rho)= \emptyset$ -- recall that
$\Omega$ is compact and so is bounded.

The set $\Bad^*(\mathcal{R}, \beta, \rho)$ is easily seen to be a
generalization of the classical set $\Bad$ of badly approximable
numbers. Recall, a real number $x$ is said to be badly approximable if
there exists a constant $c(x) > 0 $ such that $|x -p/q| \geq c(x)/q^2
$ for all rational $p/q$. A result of Jarn\'{\i}k \cite{Ja} states
that the Hausdorff dimension of $\Bad$ is maximal; i.e. $\dim \Bad
=1$.  Our initial aim is to find a suitably general framework which
allows us to conclude that $\dim \Bad^*(\mathcal{R}, \beta, \rho) =
\dim \Omega$; that is to say that the set of badly approximable points
in $\Omega$ is of maximal dimension. To a certain extent, this paper
complements \cite{BDV} in which a general framework for establishing
measure theoretic laws for `well approximable' sets is established.

A few words about our chosen notation are in order. In the above setup
and its generalization in \S\ref{gensetup}, the sets of badly
approximable elements will be denoted by $\Bad^*$ followed by the
appropriate variables in brackets. In applications we define a set,
usually denoted by $\Bad$ with appropriate arguments, and show that
this set may be realized as a specialization of a general set
$\Bad^*$.

\subsection{The conditions on the setup \label{condsetup}}

Throughout, a ball $B(c,r)$ with centre $c$ and radius $r$ is
defined to be the set $\{x \in X : d(c,x) \leq r \}$. Thus all
balls will be assumed to be closed unless stated otherwise and by
definition a ball is a subset of $X$. The following conditions on
the measure $m$ and the function $\rho$ will play a central role
in our work.

\begin{enumerate}[(A)]
\item There exist strictly positive constants $\delta$ and
$r_0$ such that for $c \in \Omega$ and $r \leq r_0$
\begin{equation*}
  a \, r^\delta \leq \; \; m(B(c,r)) \;\; \leq b \, r^\delta  \ ,
\end{equation*}
where $0 < a \leq 1 \leq b$ are constants independent
of the ball.
\end{enumerate}
It is easily verified that if the measure $m$ supported on $\Omega$ is
of type (A) then $\dim \Omega = \delta$. Trivially, this implies that
$\dim X \geq \delta$. See \S\ref{prel} for the details.

\begin{enumerate}[(B)]
\item For $ k >1 $ sufficiently large and any integer $n \geq 1$,
  \begin{equation*}
    \lambda^l (k) \leq \dfrac{\rho(k^n)}{\rho(k^{n+1})} \leq
    \lambda^u(k) \label{regk}
  \end{equation*}
  where $\lambda^l$ and $\lambda^u$ are lower and upper bounds
  depending only on $k$ such that $\lambda^l(k) \to \infty $ as $k \to
  \infty $.
\end{enumerate}
Note that this condition on $\rho$ is satisfied by any function
satisfying the following `regularity' condition. There exist a
constant $k > 1$ such that for $r$ sufficiently large
\begin{equation*}
  \lambda^l \ \leq \ \dfrac{\rho(r)}{\rho(k \, r)} \ \leq \ \lambda^u \
  \ ,
\end{equation*}
where $1< \lambda^l \leq \lambda^u $ are constants independent of $r$
but may depend on $k$.

\subsection{The  result}

First some useful notation.  For any $k > 1 $ let $B_n := \{x \in
\Omega : d(c,x)\leq \rho(k^n)\} $ denote a generic closed ball of
radius $\rho(k^n)$ with centre $c$ in $\Omega$ and for $\theta \in
\mathbb{R}^+$, let $\theta B_n := \{x \in \Omega : d(c,x)\leq \theta
\rho(k^n)\} $ denote the ball $B_n$ scaled by $\theta$. Notice, that
by definition any generic ball $B_n$ is a subset of $\Omega$.  Also,
for $n\geq 1$ let $J(n):=\{\alpha \in J: k^{n-1} \leq \beta_\alpha <
k^n\}$.

\newpage

\begin{thm}
  Let $(X,d)$ be a metric space and $(\Omega, d,m)$ a compact measure
  subspace of $X$. Let the measure $m$ and the function $\rho$ satisfy
  conditions (A) and (B) respectively.  For $k\geq k_0 >1$, suppose
  there exists some $\theta \in \mathbb{R}^+$ so that for $n\geq 1$
  and any ball $B_n$ there exists a collection $\cthbn$ of disjoint
  balls $2\theta B_{n+1}$ contained within $\theta B_n $ satisfying
  \begin{equation}
    \#\cthbn \ \geq \ \kappa_1
    \left(\dfrac{\rho(k^n)}{\rho(k^{n+1})}\right)^\delta \label{h1}
  \end{equation}
  and
  \begin{equation}
    \#   \left\{2\theta B_{n+1} \subset \cthbn : \min_{\alpha \in
        J(n+1)}d(c,R_\alpha) \leq 2\theta \rho(k^{n+1}) \right\} \ \leq \
    \kappa_2 \left(\dfrac{\rho(k^n)}{\rho(k^{n+1})}\right)^\delta  \, ,
    \label{h2}
  \end{equation}
  where $0<\kappa_2 < \kappa_1 $ are absolute constants independent of
  $k$ and $n$.  Furthermore, suppose $ \dim (\cup_{\alpha \in J}
  R_\alpha) < \delta $. Then
  \[\dim\Bad^*(\mathcal{R}, \beta, \rho)=\delta \   \ . \] \label{main}
\end{thm}

\noindent\textbf{Remarks:}
\begin{enumerate}[(i)]
\item In applications, the `scaling factor' $\theta$ is usually
  dependent on $k$ -- see the basic example below.  For $k$
  sufficiently large, it is always possible to find the collection
  $\cthbn$ satisfying condition (\ref{h1}) -- see \S\ref{prel} for the
  details. Finally, note that in the case that the resonant sets are
  points $ \dim (\cup_{\alpha \in J} R_\alpha) =0 $ and the hypothesis
  that $\dim (\cup_{\alpha \in J} R_\alpha) < \delta $ is trivially
  satisfied. This follows from the fact that the indexing set $J$ is
  countable.

\vspace*{2ex}

\item We suspect that Theorem \ref{main} can be established using
  Schmidt games \cite{MR0195595} -- a standard mechanism in the
  subject to prove such full dimension results. However, we will
  deduce the result from a more general one (Theorem \ref{maingen}
  below) which we have not been able to prove using Schmidt games. In
  fact, it is not at all clear that the Schmidt games mechanism is
  even applicable.
\vspace*{2ex}

\item Here and in subsequent theorems, we consider families of general
  resonant sets $R_\alpha$. However, in all the applications
  considered in \S5, the resonant sets are assumed to be points. There
  are natural problems of the same type where this is not the case.
  For example, when considering the classical problem of approximation
  of systems of linear forms over $\mathbb{R}$ the resonant sets are
  affine spaces in $\mathbb{R}^n$ (see \cite{MR0248090}). For reasons
  which will be explained in the final part of \S2.3, our results are
  not immediately applicable to this situation. In a forthcoming paper
  \cite{KTV}, we will treat this aspect and related problems.
\end{enumerate}

\subsection{The basic example:  $\Bad$}
\label{basiceg}

Let $I=[0,1]$ and consider the set
\begin{equation*}
  \Bad_I := \{x \in [0,1] :|x -
  p/q| > c(x)/q^2 \; {\rm for\ all\ rationals\ } \; p/q \ \ (q>0) \} \ .
\end{equation*}
This is the classical set $\Bad$ of badly approximable numbers
restricted to the unit interval. Clearly, it can be expressed in the
form $\Bad^*(\mathcal{R}, \beta, \rho)$ with $\rho(r) := r^{-2}$ and
\begin{eqnarray*}
  & &X=\Omega:= [0,1] \  , \ \  J:= \{ (p,q) \in \mathbb{N} \times
  \mathbb{N} \backslash \{0\}: p \leq q \} \ , \ \
  \\ & & \\
  & & \alpha := (p,q) \in J \ , \ \ \  \beta_{\alpha} := q \ , \ \ \
  R_{\alpha}:= p/q   \ .
\end{eqnarray*}
The metric $d$ is of course the standard Euclidean metric; $d(x,y)
:= |x-y| \, $.  Thus in this basic example, the resonant sets
$R_{\alpha}$ are simply rational points $p/q$ and the function
$\rho$ clearly satisfies condition (B). With reference to our
framework, let the measure $m$ be one--dimensional Lebesgue
measure on $I$. Thus, $\delta =1$ and $m$  clearly satisfies
condition  (A).

We show that the conditions of Theorem \ref{main} are satisfied for
this basic example. The existence of the collection $\cthbn$, where
$B_n$ is an arbitrary closed interval of length $2 \, k^{-2n}$ follows
immediately from the following simple observation. For any two
distinct rationals $p/q$ and $p'/q'$ with $k^n \leq q, q' < k^{n+1} $
we have that
\begin{equation*}
  \left|\dfrac{p}{q} - \dfrac{p'}{q'} \right| \geq \dfrac{1}{qq'} >
  k^{-2n-2} \ .
\end{equation*}
Thus, any interval $\theta B_n $ with $\theta := \tfrac{1}{2} k^{-2}$
contains at most one rational $p/q$ with $k^n \leq q < k^{n+1} $. Let
$\cthbn$ denote the collection of intervals $2 \theta B_{n+1} $
obtained by subdividing $\theta B_n $ into intervals of length $2
k^{-2n-4} $ starting from the left hand side of $\theta B_n $. Clearly
\begin{equation*}
  \# \cthbn \ \geq \ [k^2/2] \ > \ k^2/4 \ = \text{ r.h.s. of
    \eqref{h1} with \ } \kappa_1 := 1/4 \ .
\end{equation*}
Also, in view of the above observation, for $k$ sufficiently large
\begin{equation*}
  \text{l.h.s. of \eqref{h2}} \ \leq \ 1 \ < \ k^2/8 \ =
  \text{r.h.s. of \eqref{h2} with } \kappa_2 := 1/8 \ .
\end{equation*}

The upshot of this is that Theorem \ref{main} implies that
\begin{equation*}
  \dim  \ \Bad_I = \  1  \ \ .
\end{equation*}
In turn, since $\Bad$ is a subset of $\mathbb{R}$, this implies that
$\dim \Bad = 1 $ -- the classical result of Jarn\'{\i}k \cite{Ja}.

\section{A more general framework}
\label{gensetup}

We now consider a more general framework in which the `badly
approximable' set consists of points avoiding `rectangular'
neighborhoods of resonant sets rather than simply `balls'.

Let $(X,d)$ be the product space of $t$ metric spaces $(X_i,d_i)$ and
let $(\Omega, d)$ be a compact subspace of $X$ which contains the
support of a non-atomic finite measure $m$. As before, let
$\mathcal{R}=\{R_\alpha \in X : \alpha \in J\}$ be a family of subsets
$R_\alpha$ of $X$ indexed by an infinite, countable set $J$. Thus,
each resonant set $R_\alpha$ can be split into its $t$ components
$R_{\alpha,i} \subset (X_i,d_i)$. As before, let $\beta: J \rightarrow
\mathbb{R}^+ : \alpha \rightarrow \beta_\alpha$ be a positive function
on $J$ and assume that the number of $\alpha \in J$ with
$\beta_\alpha$ bounded above is finite.

For each $1 \leq i \leq t $, let $\rho_i : \mathbb{R}^+ \rightarrow
\mathbb{R}^+ : r \to \rho_i (r) $ be a real, positive function such
that $\rho_i(r) \rightarrow 0$ as $r \rightarrow \infty$ and that
$\rho_i$ is decreasing for $r$ large enough. Furthermore, assume that
$\rho_1(r) \geq \rho_2(r) \geq \cdots \geq \rho_t(r)$ for $r$ large --
the ordering is irrelevant. Given a resonant set $R_\alpha$, let
\begin{equation*}
  F_{\alpha}(\rho_1, \ldots, \rho_t) := \{x \in X \ : \ d_i(x_i,
  R_{\alpha,i}) \leq \rho_i(\beta_\alpha) \text{ for all } 1
  \leq i \leq t \} \ ,
\end{equation*}
denote the `rectangular' $(\rho_1, \ldots, \rho_t)$-neighborhood of
$R_\alpha$ and consider the set
\begin{multline*}
  \Bad^*(\mathcal{R}, \beta, \rho_1, \ldots, \rho_t) :=\{x \in \Omega
  :\exists \; c(x)>0 \text{ such that} \\
  x \notin c(x) \, F_{\alpha}(\rho_1, \ldots, \rho_t) \text{ for all }
  \alpha\in J\} \ .
\end{multline*}
Thus, $x \in \Bad^*(\mathcal{R}, \beta, \rho_1, \ldots, \rho_t)$ if
there exists a constant $c(x) > 0$ such that for all $\alpha \in J$,
\begin{equation*}
  d_i(x_i,R_{\alpha,i}) \ \geq \ c(k) \, \rho_i(\beta_\alpha)
  \quad (1\leq i\leq t) \ .
\end{equation*}

Clearly, $\Bad^*(\mathcal{R}, \beta, \rho_1,\ldots \rho_t)$ is
precisely the set $\Bad^*(\mathcal{R}, \beta, \rho)$ of \S\ref{setup}
in the case $t=1$.  The overall aim of this section is to find a
suitably general framework which gives a lower bound for the Hausdorff
dimension of $\Bad^*(\mathcal{R}, \beta, \rho_1, \ldots, \rho_t) $.
We shall assume that $\sup_{\alpha \in J} \rho_i(\beta_\alpha)$ is
finite for each $i$ without loss of generality -- otherwise
$\Bad^*(\mathcal{R}, \beta, \rho_1, \ldots, \rho_t) = \emptyset$ and
there is nothing to prove.

\subsection{The conditions on the  general framework}
\label{condgensetup}

Given $l_1, \ldots,l_t \in \mathbb{R}^+$ and $c \in \Omega $ let
\begin{equation*}
  F(c;l_1,\ldots, l_t) \ := \ \{x \in X \ : \ d_i(x_i, c_i) \leq l_i
  \text{ for all } \ 1 \leq i \leq t \} \ ,
\end{equation*}
denote the closed `rectangle' centred at $c$ with `sidelengths'
determined by $l_1,\ldots, l_t$. Also, for any $k>1$ and $n \in
\mathbb{N}$, let $F_n$ denote a generic rectangle $F(c;\rho_1(k^n),
\ldots,\rho_t(k^n) ) \cap \Omega $ in $\Omega$ centred at a point $c$
in $\Omega$. As before, $B(c,r)$ is a closed ball with centre $c$ and
radius $r$.  The following conditions on the measure $m$ and the
functions $\rho_i$ will play a central role in our general framework.
The first two are reminiscent of conditions (A) and (B) of
\S\ref{condsetup}.

\begin{enumerate}[(A*)]
\item There exists a strictly positive constant $\delta$ such that for
  any $c \in \Omega$
  \begin{equation*}
    \liminf_{r\to 0} \dfrac{\log m(B(c,r))}{\log r} \ = \ \delta \ .
  \end{equation*}
\end{enumerate}
It is easily verified that if the measure $m$ supported on $\Omega$ is
of type (A*) then $\dim \Omega \geq \delta$ \cite[Proposition
4.9]{falc} and so $\dim X \geq \delta$. Clearly condition (A) of
\S\ref{condsetup} implies (A*).
\begin{enumerate}[(B*)]
\item For $ k >1 $ sufficiently large, any integer $n \geq 1$ and any
  $i \in \{1, \ldots, t\}$,
  \begin{equation*}
    \lambda^l_i (k) \leq \dfrac{\rho_i(k^n)}{\rho_i(k^{n+1})} \leq
    \lambda^u_i(k),
    \label{regkg}
  \end{equation*}
  where $\lambda^l_i$ and $\lambda^u_i$ are lower and upper constants
  such that $\lambda^l_i(k) \to \infty $ as $k \to \infty$.
\end{enumerate}
Clearly, this is just condition (B) of \S\ref{condsetup} imposed on
each function $\rho_i$.
\begin{enumerate}[(C*)]
\item There exist constants $0<a\leq 1 \leq b $ and $l_0>0$ such that
  \begin{equation*}
    a \ \leq \ \dfrac{m(F(c;l_1, \ldots, l_t))}{m(F(c';l_1, \ldots,
      l_t))} \ \leq \ b ,
    \label{mcomp}
  \end{equation*}
  for any $c,c' \in \Omega$ and any $l_1, \ldots, l_t \leq l_0$.
\end{enumerate}
This condition implies that rectangles of the same size centred at
points of $\Omega$ have comparable $m$ measure.
\begin{enumerate}[(D*)]
\item There exist strictly positive constants $D$ and $l_0$ such that
  \begin{equation*}
    \dfrac{m(2\,F(c;l_1, \ldots, l_t))}{m(F(c;l_1, \ldots, l_t))} \
    \leq \ D,
    \label{mdoubling}
  \end{equation*}
  for any $c \in \Omega$ and any $l_1, \ldots, l_t \leq l_0$.
\end{enumerate}
This condition simply says that the measure $m$ is `doubling' with
respect to rectangles. In terms of achieving our aim of obtaining a
lower bound for $\dim \Bad^*(\mathcal{R}, \beta, \rho_1,\ldots \rho_t)
$, the above four conditions are rather natural.  The following final
condition is in some sense the only genuine technical condition and is
not particularly restrictive.
\begin{enumerate}[(E*)]
\item For $k > 1 $ sufficiently large and any integer $n \geq 1 $
  \begin{equation*}
    \dfrac{m(F_n)}{m(F_{n+1})} \ \geq \ \lambda(k) \  ,
    \label{growcond}
  \end{equation*}
  where $\lambda$ is a constant such that $\lambda(k) \to \infty $ as
  $k \to \infty$.
\end{enumerate}

\subsection{The general result}
\label{maingensetup}

Recall, that $F_n:= \{x \in \Omega : \ d_i(x_i, c_i) \leq \rho_i(k^n)
\ \ \text{for all} \ 1 \leq i \leq t \} $ is a generic rectangle with
centre $c$ in $\Omega$ and `sidelengths' determined by $ \rho_i(k^n) $
and for $\theta \in \mathbb{R}^+$, $\theta F_n $ is the rectangle $F_n
$ scaled by $\theta$. Also, for $n\geq 1$ let $J(n):=\{\alpha \in J:
k^{n-1} \leq \beta_\alpha < k^n\}$.

\begin{thm}
  Let $(X,d)$ be the product space of the metric spaces
  $(X_1,d_1),\allowbreak \ldots, (X_t,d_t)$ and let $(\Omega, d,m)$ be
  a compact measure subspace of $X$.  Let the measure $m$ and the
  functions $\rho_i$ satisfy conditions (A*) to (E*). For $k\geq k_0
  >1$, suppose there exists some $\theta \in \mathbb{R}^+$ so that for
  $n\geq 1$ and any rectangle $F_n$ there exists a disjoint collection
  $\cthfn$ of rectangles $2\theta F_{n+1}$ contained within $\theta
  F_n $ satisfying
  \begin{equation}
    \#\cthfn \ \geq \ \kappa_1 \ \dfrac{m(\theta F_n)}{m(\theta
      F_{n+1})} \label{hg1}
  \end{equation}
  and
  \begin{multline}
    \# \left\{2\theta F_{n+1} \subset \cthfn : \min_{\alpha \in
        J(n+1)}d_i(c_i,R_{\alpha,i}) \leq 2\theta \rho_i(k^{n+1}) \ \
      \text{for any} \  1 \leq i \leq t \right\}  \\
    \leq \kappa_2 \ \dfrac{m(\theta F_n)}{m(\theta F_{n+1})} \ .
    \label{hg2}
  \end{multline}
  where $0<\kappa_2 < \kappa_1 $ are absolute constants independent of
  $k$ and $n$. Furthermore, suppose $ \dim (\cup_{\alpha \in J}
  R_\alpha) < \delta $. Then
  \begin{equation*}
    \dim\Bad^*(\mathcal{R}, \beta, \rho_1,\ldots \rho_t) \geq \delta \
    \ .
  \end{equation*}
 \label{maingen}
\end{thm}

\noindent\textbf{Remarks:} For $k$ sufficiently large, it is always
possible to find the collection $\cthfn$ satisfying condition
\eqref{hg1}.  Clearly, the lower bound for $\dim \Bad^*(\mathcal{R},
\beta, \rho)$ in Theorem \ref{main} is an immediate consequence of
Theorem \ref{maingen}. To see this, simply note that if $t=1$ then the
rectangles $F_n $ are balls $B_n$ and if conditions (A) and (B) are
satisfied then trivially so are the conditions (A*) to (E*). In fact,
if condition (A*) is replaced by the stronger condition (A) in the
above theorem, then we are able to conclude that
$\dim\Bad(\mathcal{R}, \beta, \rho_1,\ldots \rho_t)= \delta$ -- see
below.

We now consider an extremely useful specialization of the above
general framework in which the space $\Omega$ is a product space
equipped with a product measure.

\begin{thm}
  For $1\leq i \leq t $, let $(X_i,d_i)$ be a metric space and
  $(\Omega_i, d_i,m_i)$ be a compact measure subspace of $X_i$ where
  the measure $m_i$ satisfies condition (A) with exponent $\delta_i$.
  Let $(X,d)$ be the product space of the spaces $(X_i,d_i)$ and let
  $(\Omega, d,m)$ be the product measure space of the measure spaces
  $(\Omega_i, d_i,m_i)$.  Let the functions $\rho_i$ satisfy condition
  (B*). For $k\geq k_0 >1$, suppose there exists some $\theta \in
  \mathbb{R}^+$ so that for $n\geq 1$ and any rectangle $F_n$ there
  exists a disjoint collection $\cthfn$ of rectangles $2\theta
  F_{n+1}$ contained within $\theta F_n $ satisfying
  \begin{equation}
    \#\cthfn \ \geq \ \kappa_1 \ \prod_{i=1}^{t}
    \left(\dfrac{\rho_i(k^n)}{\rho_i(k^{n+1})}\right)^{\delta_i}
    \label{hgs1}
  \end{equation}
  and
  \begin{multline}
    \# \left\{2\theta F_{n+1} \subset \cthfn : \min_{\alpha \in
        J(n+1)}d_i(c_i,R_{\alpha,i}) \leq 2\theta \rho_i(k^{n+1}) \ \
      \text{for any} \  1 \leq i \leq t \right\}  \\
    \leq \kappa_2 \ \prod_{i=1}^{t}
    \left(\dfrac{\rho_i(k^n)}{\rho_i(k^{n+1})}\right)^{\delta_i} \ ,
    \label{hgs2}
  \end{multline}
  where $0<\kappa_2 < \kappa_1 $ are absolute constants independent of
  $k$ and $n$. Furthermore, suppose $ \dim (\cup_{\alpha \in J}
  R_\alpha ) < \sum_{i=1}^t \delta_i $. Then
  \begin{equation*}
    \dim\Bad^*(\mathcal{R}, \beta, \rho_1,\ldots \rho_t)  =
    \sum_{i=1}^t \delta_i \   \ .
  \end{equation*}
  \label{maingencorr}
\end{thm}

The deduction of Theorem \ref{maingencorr} from Theorem \ref{maingen}
is relatively straightforward and hinges on the following simple
observation. Since $m$ is the product measure of the measures $m_i$
and the latter satisfy condition (A) with exponents $\delta_i$ ($1\leq
i \leq t$), we have for any $c \in \Omega$ and any $l_1, \ldots, t_t
\leq l_0$ that
\begin{equation}
  a^t \ \leq \ \dfrac{m(F(c;l_1, \ldots, l_t))}{\prod_{i=1}^{t}
    l_i^{\delta_i}} \ \leq \ b^t.
  \label{noworries}
\end{equation}
It follows that conditions (C*) and (D*) are trivially satisfied as is
condition (A) with $\delta:= \sum_{i=1}^t \delta_i $.  Recall, that
(A) implies (A*).  Also, (\ref{noworries}) together with (B*) implies
that condition (E*) is satisfied.  Thus, Theorem \ref{maingen} implies
the desired lower bound estimate for $ \dim\Bad^*(\mathcal{R}, \beta,
\rho_1,\ldots \rho_t)$. The complementary upper bound estimate is a
simple consequence of the fact that $m$ satisfies (A). If $m$
satisfies (A), then $\dim \Omega = \delta $ \cite[Proposition
4.9]{falc} and since $\Bad^*(\mathcal{R}, \beta, \rho_1,\ldots,
\rho_t)\subseteq \Omega$ the upper bound follows.

\subsection{The general basic example: $\Bad(i,j)$ \label{basicgen}}

For $i,j \geq 0$ with $i+j =1$, denote by $\Bad(i,j)$ the set of
$(i,j)$-badly approximable pairs $(x_1,x_2) \in \mathbb{R}^2$; that
is $(x_1,x_2) \in \Bad(i,j) $ if there exists a positive constant
$c(x_1,x_2)$ such that for all $q \in \mathbb{N}$
\begin{equation*}
  \max \{ \; ||qx_1||^{1/i} \; , \ ||qx_2||^{1/j} \,  \} \ > \
  c(x_1,x_2) \ q^{-1} \ \ \ ,
\end{equation*}
where $ || \, \cdot \, || $ denotes the distance of a real number to
the nearest integer. In the case $i=j=1/2$, the set under
consideration is simply the standard set of badly approximable pairs.
If $i=0$ we identify the set $\Bad(0,1)$ with $\mathbb{R} \times \Bad
$ where $\Bad$ is the set of badly approximable numbers. That is,
$\Bad(0,1)$ consists of pairs $(x_1,x_2)$ with $x_1 \in \mathbb{R}$
and $x_2 \in \Bad$.  The roles of $x_1$ and $x_2$ are reversed if
$j=0$.  Recently \cite{PV}, it has been shown that $\dim \Bad(i,j)=2$.
We now show that this result is in fact a simple consequence of
Theorem \ref{maingencorr}.

Let $\Bad_{I^2}(i,j) := \Bad(i,j) \cap I^2 $ where $I^2 := [0,1]
\times [0,1]$. Without loss of generality assume that $i \leq j$.
Clearly, it can be expressed in the form $\Bad^*(\mathcal{R}, \beta,
\rho_1,\rho_2)$ with $\rho_1(r) := r^{-(1+i)}$, $\rho_2(r) :=
r^{-(1+j)}$ and
\begin{eqnarray*}
  & &X=\Omega:= I^2 \  , \ \  J:= \{ ((p_1,p_2),q) \in \mathbb{N}^2
  \times \mathbb{N} \backslash \{0\}: p_1,p_2 \leq q \} \ , \ \
  \\ & & \\
  & & \alpha := ((p_1,p_2),q) \in J \ , \ \ \beta_{\alpha} := q \ , \ \
  R_{\alpha}:= (p_1/q,p_2/q)   \ .
\end{eqnarray*}
Furthermore, $d_1=d_2 $ is the standard Euclidean metric on $I$ and
$m_1=m_2$ is one--dimensional Lebesgue measure on $I$. By definition,
the metric $d$ on $I^2$ is the product metric $d_1 \times d_1 $ and
the measure $m:= m_1 \times m_1$ is simply two--dimensional Lebesgue
measure on $I^2 $.

We show that the conditions of Theorem \ref{maingencorr} are satisfied
for this basic example. Clearly the functions $\rho_1, \rho_2 $
satisfy condition (B*) and the measures $m_1,m_2$ satisfy condition
(A) with $\delta_1=\delta_2 = 1 $. We now need to establish the
existence of the collection $\cthfn$, where $F_n $ is an arbitrary
closed rectangle of size $2 k^{-n(1+i)} \times 2 k^{-n(1+j)} $.  To
start with, note that $m(\theta F_n ) = 4 \theta^2 k^{-3n}$. Now
assume there are at least three rational points $(p_1/q,p_2/q),
(p'_1/q',p'_2/q') $ and $(p''_1/q'',p''_2/q'') $ with
\begin{equation*}
  k^n \leq q,q',q'' < k^{n+1}
\end{equation*}
lying within $\theta F_n $. Suppose for the moment that they do not
lie on a line and form the triangle $ \Delta$ sub-tended by them.
Twice the area of the triangle $ \Delta$ is equal to the absolute
value of the determinant
\begin{equation*}
  \det \ := \
  \begin{vmatrix}
    1 & p_1/q & p_2/q \\
    1 & p'_1/q' & p'_2/q' \\
    1 & p''_1/q'' & p''_2/q''
  \end{vmatrix}   \ .
\end{equation*}

Then, in view of the denominator constraint, it follows that
\begin{equation*}
  2 \,\times m(\Delta) \geq \dfrac{1}{q q' q''} > k^{-3(n+1)} \ .
\end{equation*}
Now put
\begin{equation*}
  \theta := 2^{-1} (2 k^3)^{-1/2}.
\end{equation*}
Then $m(\Delta) > m( \theta F_n)$ and this is impossible since $\Delta
\subset \theta F_n$.  The upshot of this is that the triangle in
question can not exist. Thus, if there are two or more rational points
with $ k^n \leq q < k^{n+1} $ lying within $\theta F_n$ then they must
lie on a line $\mathcal{L}$.

Starting from a `corner' of the rectangle $\theta F_n$, partition
$\theta F_n$ into rectangles $2\theta F_{n+1} $ of size $4
k^{-(n+1)(1+i)} \times 4 k^{-(n+1)(1+j)} $ and denote by $\cthfn$ the
collection of rectangles $2\theta F_{n+1}$ obtained. Trivially
\begin{equation*}
  \# \cthfn \ \geq \ \left[\dfrac{2 \theta k^{-n(1+i)} }{4 \theta
      k^{-(n+1)(1+i)} }\right] \ \left[\dfrac{2 \theta k^{-n(1+j)} }{4
      \theta k^{-(n+1)(1+j)} } \right] \ \geq \dfrac{k^3}{16} \ .
\end{equation*}
In view of the above `triangle' argument we have that
\begin{multline*}
  \#  \left\{2\theta F_{n+1} \subset \cthfn : \min_{\alpha \in
      J(n+1)}d_i(c_i,R_{\alpha,i})  \leq   2\theta \rho_i(k^{n+1}) \  \
    \text{for all} \  1 \leq i \leq t \right\} \\ \\ \vspace*{-5ex}
  \leq  \ \# \left\{2\theta F_{n+1} \subset \cthfn : 2\theta F_{n+1} \cap
    \mathcal{L} \neq \emptyset \right\} \ ,
\end{multline*}
where $\mathcal{L}$ is any line passing through $\theta F_n$. Recall,
that we are assuming that $i \leq j$. A simple geometric argument
ensures that for $k$ sufficiently large
\begin{eqnarray*}
  \# \left\{2\theta F_{n+1} \subset \cthfn : 2\theta F_{n+1} \cap
    \mathcal{L} \neq \emptyset \right\}   \ & \leq & \   \left[\dfrac{2
      \theta k^{-n(1+j)} }{4 \theta k^{-(n+1)(1+j)} } \right] =
  \left[\dfrac{ k^{1+j}}{2} \right] \\ & & \\  \ & \leq & \   k^{1+j} \ \leq \
  k^3/32 \ .
\end{eqnarray*}

The upshot of this is that the collection $\cthfn $ satisfies the
required conditions and Theorem \ref{maingencorr} implies that
\begin{equation*}
  \dim \Bad_{I^2}(i,j) \ = \ 2 \ .
\end{equation*}
In turn, since $\Bad(i,j)$ is a subset of $\mathbb{R}^2$, this implies
that $\dim \Bad(i,j) = 2 $.

In \cite{PV}, the stronger result that $\dim \Bad(i,j) \cap \Bad(1,0)
\cap \Bad(0,1) = 2 $ is established; i.e. the set of pairs $(x_1,x_2)$
with $x_1$ and $x_2$ both badly approximable numbers and an
$(i,j)$-badly approximable pair has full dimension. In
\S\ref{badijo}, we obtain a much more general result and remark on a
beautiful conjecture of W.M. Schmidt. In full generality, Schmidt's
conjecture states that $ \Bad(i,j) \cap \Bad(i',j') \neq \emptyset $.
It is a simple exercise to show that if Schmidt's conjecture is false
for some pairs $(i,j)$ and $(i',j')$ then Littlewood's conjecture in
simultaneous Diophantine approximation is true.

We now turn our attention to the natural generalization of $\Bad(i,j)$
to higher dimensions.  For any $N$-tuple of real numbers $i_1,
\ldots,i_N \geq 0 $ such that $\sum i_r = 1 $, denote by $\Bad(i_1,
\ldots,i_N)$ the set of points $(x_1, \ldots,x_N) \in \mathbb{R}^N$
for which there exists a positive constant $ c(x_1, \ldots, x_N)$ such
that for any $q \in \mathbb{N}$,
\begin{equation*}
  \max \{ \; ||qx_1||^{1/i_1} \; , \ldots, \ ||qx_N||^{1/i_N} \, \} \
  > \ c(x_1, \ldots,x_N) \ q^{-1}.
\end{equation*}
Clearly, the two-dimensional argument can easily be modified to show
that
\begin{equation*}
  \dim \Bad(i_1, \ldots,i_N) = N \ .
\end{equation*}
The key modification is the following lemma which naturally extends
the main feature of the `triangle' argument in dimension two to a
`simplex' one in dimension $N$.

\begin{lem}{(Simplex Lemma)}
  Let $N \geq 1$ be an integer and $k > 1$ be a real number.  Let $E
  \subseteq \mathbb{R}^N$ be a convex set of $N$-dimensional Lebesgue
  measure
  \begin{equation*}
    \vert E \vert \ \leq \ (N! \, )^{-1} k^{-(N+1)} \ .
  \end{equation*}
  Suppose that $E$ contains $N+1$ rational points $(p_i^{(1)}/q_i,
  \ldots, p_i^{(N)}/q_i)$ with $1 \leq q_i < k$, where $0 \leq i \leq
  N$. Then these rational points lie in some hyperplane.
\end{lem}

\begin{pf*}{\textbf{Proof.} \ }
  Suppose to the contrary that this is not the case. In that case, the
  rational points $(p_i^{(1)}/q_i, \ldots, p_i^{(N)}/q_i)$ where $0
  \leq i \leq N$ are distinct.  Consider the $N$-dimensional simplex
  $\Delta $ subtended by them; i.e. an interval when $N=1$, a triangle
  when $N=2$, a tetrahedron when $N=3$ and so on. Clearly, $\Delta$ is
  a subset of $E$ since $E$ is convex.  The volume of the simplex
  $\vert \Delta \vert$ times $N$ factorial is equal to the absolute
  value of the determinant
  \begin{equation*}
    \det :=
    \begin{vmatrix}
      1 & p_0^{(1)}/q_0 & \cdots & p_0^{(N)}/q_0 \\
      1 & p_1^{(1)}/q_1 & \cdots & p_1^{(N)}/q_1 \\
      \vdots & \vdots & & \vdots \\
      1 & p_N^{(1)}/q_N & \cdots & p_N^{(N)}/q_N \\
    \end{vmatrix}.
  \end{equation*}
  As this determinant is not zero, it follows from the assumption made
  on the $q_i$ that
  \begin{equation*}
    N! \times \vert \Delta \vert = \vert \det \vert \geq \dfrac{1}{q_0
      q_1 \cdots q_N} > k^{-(N+1)} \, .
  \end{equation*}
  Consequently, $\vert \Delta \vert > (N! \, )^{-1} k^{-(N+1)} \geq
  \vert E \vert $. This contradicts the fact that $\Delta \subseteq
  E$.
\hfill $ \Box $
\end{pf*}

\noindent\textbf{Remarks:}
\begin{enumerate}[(i)]
\item The Simplex Lemma should be viewed as the higher dimensional
  generalization of the following simple fact already exploited in the
  \S\ref{basiceg}: on the real line $\mathbb{R}$ an interval $I_k$ of
  length $1/k^2$ can contain at most one rational $p/q$ with $1\leq q
  < k $.  This follows from the trivial observation that if $1 \leq q,
  q' < k$ then $|p/q-p'/q'| \geq 1/qq' > 1/ k^2$.
\vspace*{2ex}

\item Our general setup will be applied to settings other than subsets
  of $\mathbb{R}^N$ (see \S\ref{applic}). In most of these, an
  analogue of the Simplex Lemma will be required. In these settings we
  will either give a complete proof or sketch the argument required in
  two dimensions; i.e. the analogue of the `triangle' argument.  Based
  on the proof of the Simplex Lemma in $\mathbb{R}^N$, it should then be
  obvious how to extend the $N=2$ argument to higher dimensions. In
  short, within this paper the main ideas are always exposed on
  establishing a given $N$-dimensional statement in the $N=2$ case.
  The proof in higher dimensions requires no new ideas.  Thus in all
  the various applications of our general framework, for the sake of
  both clarity and notation we shall stick to $N=2$ in proofs.
\vspace*{2ex}

\item The `triangle' argument (or variants thereof) described above is
  critical in most of the applications considered in this paper (see
  \S\ref{applic}). To some extent this is the reason why our main
  results cannot be directly applied to the problem of badly
  approximable systems of linear forms.  In this case the resonant
  sets $R_\alpha$ are affine spaces and although the `triangle' or
  more generally the `simplex' approach remains the main ingredient
  it requires deeper considerations in the geometry of numbers to
  successfully execute it. We will return to this and other aspects of
  the linear forms theory in a forthcoming paper \cite{KTV}.
\end{enumerate}

\section{Preliminaries}
\label{prel}

In this short section we define Hausdorff measure and dimension in
order to establish some notation and then describe a method for
obtaining lower bounds for the dimension.

Suppose $\Omega$ is a non--empty subset of $(X,d)$. For $\rho > 0$, a
countable collection $ \left\{B_{i} \right\} $ of balls in $X$ with
radii $r_i \leq \rho $ for each $i$ such that $\Omega \subset
\bigcup_{i} B_{i} $ is called a $ \rho $-cover for $\Omega$. Clearly
such a cover always exists for totally bounded metric spaces. Let $s$
be a non-negative number and define
\begin{equation*}
  \mathcal{H}^{s}_{\rho} (\Omega) \,
  = \, \inf \left\{ \sum_{i} r_i^s \ : \{ B_{i} \} \text{ is a
      $\rho$-cover of } \Omega \right\} \, \ ,
\end{equation*}
where the infimum is over all $\rho$-covers.  The
\emph{$s$-dimensional Hausdorff measure} $\mathcal{H}^{s} (\Omega)$
of $\Omega$ is defined by
\begin{equation*}
  \mathcal{H}^{s} (\Omega) := \lim_{\rho \rightarrow 0}
  \mathcal{H}^{s}_{\rho} (\Omega) \; = \; \sup_{\rho > 0 }
  \mathcal{H}^{s}_{\rho} (\Omega) \;
\end{equation*}
and the Hausdorff dimension $\dim \Omega$ of a set $\Omega$ by
\begin{equation*}
  \dim\, \Omega \, := \, \inf \left\{ s : \mathcal{H}^{s} (\Omega) =0
  \right\} = \sup \left\{ s : \mathcal{H}^{s} (\Omega) = \infty
  \right\}.
\end{equation*}
In particular when $s$ is an integer $\mathcal{H}^s$ is comparable to
$s$-dimensional Lebesgue measure. For further details see
\cite{falc,mat}. A general and classical method for obtaining a lower
bound for the Hausdorff dimension of an arbitrary set $\Omega$ is the
following mass distribution principle (see e.g. \cite[page 55]{falc}).

\begin{lem}{(Mass Distribution Principle)}
  Let $ \mu $ be a probability measure supported on a subset $\Omega$
  of $ (X,d) $.  Suppose there are positive constants $c$ and $r_0$
  such that
  \begin{equation*}
    \mu ( B ) \leq \, c \; r^s \; ,
  \end{equation*}
  for any ball $B$ with radius $r \leq r_0 \, $. Then $\mathcal H^{s}
  (\Omega) \geq 1/c \, $. In particular, we have that $\dim \Omega \geq s
  $.
\end{lem}

The following rather simple covering result will be crucial to our
proof of Theorem  \ref{maingen}.

\begin{lem}{(Covering Lemma)}
  Let $(X,d)$ be the product space of the metric spaces $(X_1,d_1),$ $
  \ldots, (X_t,d_t)$ and $\mathcal{F}$ be a finite collection of
  `rectangles' $F:=F(c;l_1, \ldots, l_t)$ with $c \in X$ and $l_1,
  \ldots, l_t$ fixed. Then there exists a disjoint sub-collection
  $\{F_m\}$ such that
  \begin{equation*}
    \bigcup_{F\in\mathcal{F}} \, F \ \subset \ \bigcup_m \ {3F}_{\! m}
    \ .
  \end{equation*}
\end{lem}

\begin{pf*}{\textbf{Proof.}  \ }
  Let $S$ denote the set of centres $c$ of the rectangles in
  $\mathcal{F}$.  Choose $c(1) \in S$ and for $k \geq 1 $,
  \begin{equation*}
    c(k+1) \ \in \ S \ \backslash \ \bigcup_{m=1}^{k} 2\, F(c(m) ;l_1,
    \ldots, l_t)
  \end{equation*}
  as long as $S \ \backslash \ \bigcup_{m=1}^{k} 2\, F(c(m) ;l_1,
  \ldots, l_t) \neq \emptyset $.  Since $\#S$ is finite, there exists
  $k_1 \leq \#S $ such that
  \begin{equation*}
  S \ \subset \ \bigcup_{m=1}^{k_1} 2\, F(c(m) ;l_1, \ldots, l_t) \ .
  \end{equation*}
  By construction, any rectangle $ F(c;l_1, \ldots, l_t)$ in the
  original collection $\mathcal{F}$ is contained in some rectangle
  $3\, F(c(m) ;l_1, \ldots, l_t)$ and since $d_i(c_i(m),c_i(n)) >
  2l_i$ for each $1\leq i \leq t $ the chosen rectangles $ F(c(m)
  ;l_1, \ldots, l_t)$ are clearly disjoint.

\hfill $ \Box $
\end{pf*}

We end this section by making use of the covering lemma to establish
the following assertion made in \S\ref{maingensetup}.  The result is
extremely useful when it comes to applying our theorems -- see
\S\ref{applic}.  With reference to Theorem \ref{maingen}, it
guarantees the existence of a disjoint collection $\cthfn$ of
rectangles with the necessary cardinality.

\begin{lem}
  Let $(X,d)$ be the product space of the metric spaces $(X_1,d_1),
  \ldots,\allowbreak  (X_t,d_t)$ and let $(\Omega, d,m)$ be a compact measure
  subspace of $X$.  Let the measure $m$ and the functions $\rho_i$
  satisfy conditions (B*) to (D*). Let $k$ be sufficiently large. Then
  for any $\theta \in \mathbb{R}^+$ and for any rectangle $F_n$ ($n
  \geq 1$) there exists a disjoint collection $\cthfn$ of rectangles
  $2\theta F_{n+1}$ contained within $\theta F_n $ satisfying
  (\ref{hg1}) of Theorem \ref{maingen}.
  \label{crucial}
\end{lem}

\begin{pf*}{\textbf{Proof.}  \ } Begin by choosing $k$ large enough so that for any $i \in \{1,
  \ldots, t\}$,
  \begin{equation}
    \dfrac{\rho_i(k^n)}{\rho_i(k^{n+1})} \ \geq 4.
    \label{large}
  \end{equation}
  That this is possible follows from the fact that $\lambda^l_i(k) \to
  \infty $ as $k \to \infty$ (condition (B*)).  Take an arbitrary
  rectangle $ F_n$ and let $l_i(n) := \theta\rho_i(k^n)$. Thus $
  \theta F_n := F(c;l_1(n), \ldots, l_t(n))$. Consider the rectangle
  $T_n \subset \theta F_n$ where
  \begin{equation*}
    T_n := F(c;l_1(n) - 2l_1(n+1), \ldots, l_t(n) - 2l_t(n+1)) \ .
  \end{equation*}
  Note that in view of (\ref{large}) we have that $T_n \supset
  \tfrac{1}{2} \theta F_n$.  Now, cover $T_n$ by rectangles $2\theta
  F_{n+1}$ with centres in $\Omega \cap T_n$. By construction, these
  rectangles are contained in $\theta F_n$ and in view of the covering
  lemma there exists a disjoint sub-collection $\cthfn$ such that
  \begin{equation*}
    T_n \subset \bigcup\limits_{2\theta F_{n+1}\subset \cthfn} 6\theta
  F_{n+1} \ \ .
  \end{equation*}

  Using that fact that rectangles of the same size centred at points
  of $\Omega$ have comparable $m$ measure (condition (C*)), it follows
  that
  \begin{equation*}
    a \, m(\tfrac{1}{2} \theta F_n) \ \leq \ m(T_n) \leq \#\cthfn \ b
    \, m(6\theta F_{n+1}) \ \ .
  \end{equation*}
  Using that fact that the measure $m$ is doubling on rectangles
  (condition (D*)), so that
  $m(\tfrac{1}{2} \theta F_n) \geq D^{-1} m(\theta F_n)$ and
  $m(6\theta F_{n+1}) \leq m(8\theta F_{n+1}) \leq D^3m(\theta
  F_{n+1})$, it follows that
  \begin{equation*}
    \#\cthfn \ \geq \
  \dfrac{a}{bD^4} \ \ \dfrac{m(\theta F_{n})}{m(\theta F_{n+1})} \ \ .
  \end{equation*}
\hfill $ \Box $
\end{pf*}

\noindent \textbf{Remark.} Clearly, with reference to Theorem
\ref{main}, the above lemma guarantees the existence of the collection
$\cthbn$ satisfying \eqref{h1}.

\section{Proof of Theorem \ref{maingen}}

The overall strategy is as follows. For any $k$ sufficiently large we
construct a Cantor-type set $\mathbf{K}_{c(k)}$ such that
$\mathbf{K}_{c(k)}$ with at most a finite number of points removed is
a subset of $\Bad^*(\mathcal{R}, \beta, \rho_1, \ldots ,\rho_t)$.
Next, we construct a measure $\mu$ supported on $\mathbf{K}_{c(k)}$
with the property that for any ball $A$ with radius $r(A)$
sufficiently small
\begin{equation*}
\mu(A) \ \ll \ r(A)^{\delta - \epsilon(k)} \ ;
\end{equation*}
where $\epsilon(k) \to 0 $ as $ k \to \infty $.  Hence, by
construction and the mass distribution principle we have that
\begin{equation*}
  \dim \Bad^*(\mathcal{R}, \beta, \rho_1, \ldots ,\rho_t) \ \geq \
  \dim \mathbf{K}_{c(k)} \ \geq \ \delta - \epsilon(k) \ .
\end{equation*}
Now suppose that $\dim \Bad^*(\mathcal{R}, \beta, \rho_1,
\ldots,\rho_t) < \delta $.  Then, $\dim \Bad^*(\mathcal{R}, \beta,
\rho_1, \allowbreak \ldots ,\rho_t) = \delta - \eta $ for some $\eta > 0$.
However, by choosing $k$ large enough so that $\epsilon(k) < \eta $ we
obtain a contradiction and thereby the lower bound result follows.

\subsection{ The Cantor-type set $\mathbf{K}_{c(k)}$}

Choose $k_0$ sufficiently large so that for $k \geq k_0$, $\rho_i(k) $
($1\leq i \leq t$) is decreasing and the hypotheses of the theorem are
valid. Now fix some $k \geq k_0$ and suppose that
\begin{equation}
  \left\{ \alpha \in J : \beta_{\alpha} < k  \right\} \ = \
  \emptyset   \ . \label{nok}
\end{equation}
Define $\mathcal{F}_1$ to be any rectangle $\theta F_1 $ of radius
$\theta \rho(k)$ and centre $c$ in $\Omega$. The idea is to establish,
by induction on $n$, the existence of a collection $\mathcal{F}_n$ of
disjoint rectangles $\theta F_n $ such that $\mathcal{F}_n$ is nested
in $\mathcal{F}_{n-1}$; that is, each rectangle $\theta F_n $ in
$\mathcal{F}_n$ is contained in some rectangle $\theta F_{n-1} $ of
$\mathcal{F}_{n-1}$.  Also, any $\theta F_n $ in $\mathcal{F}_n$ will
have the property that for all points $x \in \theta F_n $, for all $i
\in \{1, \ldots, t\}$ and for all $\alpha \in J$ with $\beta_\alpha <
k^n$,
\begin{equation}
  d_i(x,R_{\alpha,i})  \ \geq \  c(k) \, \rho_i(\beta_\alpha),
  \label{C}
\end{equation}
where the constant
\begin{equation*}
  c(k) \ := \ \min_{1\leq i \leq t} (\theta/\lambda^u_i(k))
\end{equation*}
is dependent on $k$ but is independent of $n$.  Then, since the
rectangles $\theta F_n $ of $\mathcal{F}_n$ are closed, nested and the
space $\Omega$ is compact, any limit point in $\theta F_n $ will
satisfy \eqref{C} for all $\alpha $ in $J$ with $\beta_\alpha \geq k
$. In particular, we put
\begin{equation*}
  \mathbf{K}_{c(k)} \ := \ \bigcap\limits_{n=1}^{\infty}\mathcal{F}_n
  \ .
\end{equation*}
By construction, we have that $\mathbf{K}_{c(k)}$ is a subset of
$\Bad^*(\mathcal{R}, \beta, \rho_1,\ldots,\rho_t)$ under the
assumption (\ref{nok}).

\noindent \emph{The induction. \ } For $n=1$, (\ref{C}) is
trivially satisfied for $\mathcal{F}_1 = \theta F_1 $ since we are
assuming \eqref{nok}. Given $\mathcal{F}_n$ satisfying \eqref{C} we
wish to construct a nested collection $\mathcal{F}_{n+1}$ for which
\eqref{C} is satisfied for $n+1$. Consider any rectangle $\theta F_n
\subset \mathcal{F}_n $. We construct a `local' collection
$\mathcal{F}_{n+1}(\theta F_n) $ of disjoint rectangles $\theta
F_{n+1} $ contained in $\theta F_n$ so that for any point $x \in
\theta F_{n+1} $ the condition given by \eqref{C} is satisfied for
$n+1$.  Given that any rectangle $\theta F_{n+1} $ of
$\mathcal{F}_{n+1}(\theta F_n) $ is to be nested in $\theta F_n $, it
is enough to show that for any point $x \in \theta F_{n+1} $ the
inequalities
\begin{equation*}
  d_i(x_i,R_{\alpha,i}) \ \geq \ c(k) \, \rho_i(\beta_\alpha) \quad (1
  \leq i \leq t)
\end{equation*}
are satisfied for $\alpha \in J$ with $k^n \leq \beta_\alpha < k^{n+1}
$; i.e. with $\alpha \in J(n+1)$.

For $k$ sufficiently large, by the hypotheses of the theorem, there
exists a disjoint sub-collection $\mathcal{G}(\theta F_n)$ of
$\mathcal{C}(\theta F_n)$ of rectangles $2\theta F_{n+1} \subset
\theta F_n$ with
\begin{equation}
  \#\mathcal{G}(\theta F_n) \ = \ \left[\kappa  \ \dfrac{m(\theta
      F_n)}{m(\theta F_{n+1})}  \right] \hspace{2cm} \kappa:= \min\{1,
  \tfrac{1}{2} (\kappa_1 - \kappa_2)\} \ ,
  \label{count}
\end{equation}
and such that for any rectangle $2\theta F_{n+1} \subset
\mathcal{G}(\theta F_n) $ with centre $c$
\begin{equation*}
  \min_{\alpha \in J(n+1)} d_i(c_i,R_{\alpha,i}) \ \geq \ 2 \, \theta
  \, \rho_i(k^{n+1}) \ .
\end{equation*}
Clearly, by choosing $k$ large enough we can ensure that
$\#\mathcal{G}(\theta F_n) > 1 $ -- this makes use of conditions (D*)
and (E*).  Now let
\begin{equation*}
  \mathcal{F}_{n+1}(\theta F_n) \ := \ \left\{ \theta F_{n+1} : 2
    \theta F_{n+1} \subset \mathcal{G}(\theta F_n) \right\} \ .
\end{equation*}
Thus the rectangles of $\mathcal{F}_{n+1}(\theta F_n)$ are precisely
those of $\mathcal{G}(\theta F_n) $ but scaled by a factor $1/2$.
Then, by construction for any $ x \in \theta F_{n+1} \subset
\mathcal{F}_{n+1}(\theta F_n) $ and $1\leq i \leq t$
\begin{eqnarray*}
  d_i(x_i,R_{\alpha,i}) \ \geq \ \theta\rho_i(k^{n+1}) \ =  \
  \theta\rho_i(k^n) \ \dfrac{\rho_i(k^{n+1})}{\rho_i(k^n)}&\  \geq
  \ &
  \dfrac{\theta}{\lambda^u_i(k)}  \ \rho_i(\beta_\alpha) \\ & & \\  & \ \geq  \ & c(k) \
  \rho_i(\beta_\alpha).
\end{eqnarray*}
Here we have made use of condition (B*) and the fact that $\rho_i(k) $
is decreasing for $k \geq k_0$ and that $\alpha \in J(n+1)$. Finally
let
\begin{equation*}
  \mathcal{F}_{n+1}:=\bigcup\limits_{\theta F_n \in \mathcal{F}_n}
  \mathcal{F}_{n+1} (\theta F_n) \ .
\end{equation*}
This completes the proof of the induction step and so the construction
of the Cantor-type set
\begin{equation*}
  \mathbf{K}_{c(k)} \ := \ \bigcap_{n=1}^{\infty}\mathcal{F}_n  \ ,
\end{equation*}
where $c(k):= \min_{1\leq i \leq t} (\theta/\lambda^u_i(k)) $ and $k$
is sufficiently large.

Note, that in view of \eqref{count} we have that for $n \geq 2$
\begin{multline}
  \# \mathcal{F}_n = \#\mathcal{F}_{n-1} \ \times \
  \#\mathcal{F}_n(\theta F_{n-1}) \ = \ \prod_{m=2}^{n} \#
  \mathcal{F}_m(\theta F_{m-1}) \\
  \geq \prod_{m=2}^{n} \dfrac{\kappa}{2} \ \ \dfrac{m(\theta
    F_{m-1})}{m(\theta F_{m})} \ = \ \left( \dfrac{\kappa}{2}
  \right)^{n-1} \ \dfrac{m(\theta F_1) }{m(\theta F_{n})} \label{lbd}
  .
\end{multline}

\subsection{The measure $\mu$ on $\mathbf{K}_{c(k)}$}

We now describe a probability measure $\mu$ supported on the
Cantor--type set $\mathbf{K}_{c(k)}$ constructed in the previous
subsection. For any rectangle $\theta F_n$ in $\mathcal{F}_n$ we
attach a weight $ \mu (\theta F_n) $ which is defined recursively as
follows: for $n=1$,
\begin{equation*}
\mu(\theta F_1) \; := \; \dfrac{1}{\# \mathcal{F}_1} = 1
\end{equation*}
and for $n\geq 2$,
\begin{equation*}
  \mu(\theta F_n) \; := \; \dfrac{1}{\# \mathcal{F}_n(\theta F_{n-1})}
  \ \mu(\theta F_{n-1}) \quad (F_n \subset F_{n-1}) \; \; .
\end{equation*}
This procedure thus defines inductively a mass on any rectangle used
in the construction of $\mathbf{K}_{c(k)}$. In fact a lot more is true
-- $\mu$ can be further extended to all Borel subsets $A$ of $\Omega$
to determine $\mu(A)$ so that $\mu$ constructed as above actually
defines a measure supported on $\mathbf{K}_{c(k)}$; see
\cite[Proposition 1.7]{falc}. We state this formally as a

\noindent\textbf{Fact.}\ \ \ { \em The probability measure $\mu$
constructed above is supported on $\mathbf{K}_{c(k)}$ and for any
Borel subset $A$ of $\Omega$
\begin{equation*}
  \mu(A) \; = \; \inf\;\sum_{F\in \mathcal{F}}\mu(F) \  .
\end{equation*}
The infimum is over all coverings $\mathcal{F}$ of $A$ by
rectangles $F \in \{\mathcal{F}_n : n \geq 1\}$. }

Notice that, in view of \eqref{lbd}, we simply have that
\begin{equation*}
  \mu(\theta F_n)=\frac{1}{\#\mathcal{F}_n} \quad (n\geq 1) \; .
\end{equation*}

\subsection{A lower bound for $\dim \,\mathbf{K}_{c(k)}$}

Let $A$ be an arbitrary ball with centre $a$ not necessarily in
$\Omega$ and of radius $r(A)< \theta \rho_*(k^{n_0})$ where $\rho_*(r)
:= \max_{1\leq i \leq t}\rho_i(r) $ and $n_0$ is to be determined
later. We now determine an upper bound for $\mu(A)$ in terms of its
radius. Choose $n \geq n_0$ so that
\begin{equation*}
  \theta \rho_*(k^{n+1}) \ < \ r(A) \ \leq \ \theta \rho_*(k^n) \ .
\end{equation*}
Without loss of generality, assume that $A \cap \mathbf{K}_{c(k)} \neq
\emptyset$ since otherwise there is nothing to prove. Clearly
\begin{equation*}
  \mu(A) \ \leq \ \mathcal{N}_{n+1}(A) \times \mu(\theta F_{n+1})
\end{equation*}
where
\begin{equation*}
  \mathcal{N}_{n+1}(A) \ := \ \#\{\theta F_{n+1} \subset
  \mathcal{F}_{n+1}: \theta F_{n+1} \cap A \neq \emptyset \} \ .
\end{equation*}
If $\theta F_{n+1} \cap A \neq \emptyset $, then $ \theta F_{n+1}
\subset 3\, A $ since $ r(A) \geq \theta \rho_i(k^{n+1})$ for $1 \leq
i \leq t $. The balls in $\mathcal{F}_{n+1}$ are disjoint and have
comparable $m$ measure (condition (C*)), thus
\begin{equation*}
  \mathcal{N}_{n+1}(A) \ \leq \ \dfrac{m(3A)}{a \, m(\theta F_{n+1})}
  \ .
\end{equation*}
It follows by \eqref{lbd}, that
\begin{equation*}
  \mu(A) \ \leq \ \dfrac{m(3A)}{a \, m(\theta F_{n+1})} \times
  \dfrac{1}{\#\mathcal{F}_{n+1} } \ \ \leq \ \ \dfrac{m(3A)}{a \,
    m(\theta F_{1})} \left(\frac{2}{\kappa}\right)^{n} \ .
\end{equation*}
Using the fact that $ \rho_*(k^n) \leq \lambda^l_*(k)^{-(n-1)}
\rho_*(k) $, it is easily verified that
\begin{equation*}
  \dfrac{1}{a \, m(\theta F_{1})} \left(\frac{2}{\kappa}\right)^{n} \
  < \ \left( \dfrac{1}{ \theta \rho_*(k^n)} \right)^{\epsilon(k)}
\end{equation*}
for
\begin{equation*}
  n \ \geq \ n_1 \ := \ \left[4 + \dfrac{\log \tfrac{(\theta \,
        \rho_*(k))^{\epsilon(k)} }{a \, m(\theta F_1) } }{\log
      \frac{2}{\kappa} } \right] \quad \text{ and \ }
  \quad \epsilon(k) \ := \ \dfrac{4 \, \log \tfrac{2}{\kappa}}{\log
    \lambda^l_*(k)} \ .
\end{equation*}
Hence,
\begin{equation*}
  \mu(A) \ \leq \ m(3A) \times (\theta \rho_*(k^n) )^{- \epsilon(k)} \
  .
\end{equation*}
Since $A \cap \mathbf{K}_{c(k)} \neq \emptyset$, there exists some
point $x \in A \cap \Omega $. Moreover, $3A \subset B(x, 4 \, r(A)) $
which together with condition (A*) implies that
\begin{equation*}
  m(3A) \leq m( B(x, 4\, r(A)) ) \leq r(A)^{\delta - \epsilon(k)}
\end{equation*}
for $r(A) \leq r_0 := r_0(\epsilon(k))$. Now $\rho_*(r) \to 0 $ as $r
\to \infty$, so $\theta \rho_*(k^n) < r_0 $ for $n \geq n_2$.  Thus,
for $n \geq n_0:= \max\{n_1,n_2\}$
\begin{equation*}
  \mu(A) \ \leq \ r(A)^{\delta - \epsilon(k) } \times (\theta
  \rho_*(k^n) )^{- \epsilon(k)} \ .
\end{equation*}
On using the fact that $r(A) \leq \theta \rho_*(k^n) $, we obtain that
\begin{equation*}
  \mu(A) \ \leq \ r(A)^{\delta - 2\epsilon(k) } \ .
\end{equation*}

This together with the mass distribution principle implies that
\begin{equation*}
  \dim \mathbf{K}_{c(k)} \ \geq \ \delta - 2\epsilon(k) \ .
\end{equation*}
Note that since $\epsilon(k) \to 0 $ as $k \to \infty $ we have that
$\dim \mathbf{K}_{c(k)} \to \delta $ as $k \to \infty $.

\subsection{Completion of proof}

Recall, that $\dim(\cup_{\alpha \in J}R_\alpha) < \delta$. Now
suppose that $$\dim \Bad^*(\mathcal{R}, \beta, \rho_1, \ldots,
\allowbreak \rho_t) \  < \ \delta \ . $$ It follows that $\max\{
\dim \Bad^*(\mathcal{R}, \beta, \rho_1, \ldots,\rho_t),
\dim(\cup_{\alpha \in J} R_\alpha) \} = \delta - \eta $ for some
$\eta > 0$. Fix some $k$ sufficiently large so that $2 \,
\epsilon(k) < \eta $. Then,
\begin{equation*}
\dim \mathbf{K}_{c(k)} \ \geq \ \delta - 2 \epsilon(k) \ > \ \delta -
\eta \ .
\end{equation*}
By construction, for any point $x \in \mathbf{K}_{c(k)}$ we have for
all $\alpha \in J$ with $\beta_\alpha \geq k$ that
\begin{equation*}
   d_i(x_i,R_{\alpha,i})  \ \geq \ c(k) \, \rho_i(\beta_\alpha) \quad
   (1\leq i\leq t) \ .
\end{equation*}
Now let $J_k:= \{\alpha \in J : \beta_\alpha < k \} $.  If \eqref{nok}
is true for our fixed $k$ then $ J_k = \emptyset $ and clearly
$\mathbf{K}_{c(k)} \subseteq \Bad^*(\mathcal{R}, \beta, \rho_1,
\ldots,\rho_t) $.  In turn, $ \dim \Bad^*(\mathcal{R}, \beta, \rho_1,
\ldots,\rho_t) \geq \dim \mathbf{K}_{c(k)} > \delta - \eta $ and we
have a contradiction. So suppose, $ J_k \neq \emptyset $ and let
$\mathcal{R}_k := \{ R_\alpha: \alpha \in J_k\}$. For any fixed $k$
the number of elements in $J_k$ is finite.  So, if $x \notin R_k $
then there exists a constant $c'(x) > 0 $ such that for all $\alpha
\in J_k$,
\begin{equation*}
  d_i(x_i,R_{\alpha,i})  \ \geq \ c'(k) \, \rho_i(\beta_\alpha)
  \quad (1\leq i\leq t) \ .
\end{equation*}
Thus, for $x \in \mathbf{K}_{c(k)}\!\!\setminus \mathcal{R}_k$ and
$\alpha \in J$,
\begin{equation*}
  d_i(x_i,R_{\alpha,i})  \ \geq \ c^*(k) \, \rho_i(\beta_\alpha) \quad
  (1\leq i\leq t)\ ,
\end{equation*}
where $c^*(x) := \min\{c(k),c'(x)\}$. It follows that
$$\Bad^*(\mathcal{R}, \beta, \rho_1, \ldots,\rho_t) \supseteq
\mathbf{K}_{c(k)} \!\! \setminus \mathcal{R}_k \ , $$ and since
$\dim \mathcal{R}_k < \dim \mathbf{K}_{c(k)}$ we have that
\begin{eqnarray*}
  \dim \Bad^*(\mathcal{R}, \beta, \rho_1, \ldots,\rho_t) \ &  \geq  & \ \dim
  (\mathbf{K}_{c(k)} \!\! \setminus \mathcal{R}_k ) \\ & & \\  \  &
  =  & \  \dim \mathbf{K}_{c(k)} \  \geq  \ \delta - 2\epsilon(k) \ > \
  \delta - \eta \ .
\end{eqnarray*}
This is a contradiction and completes the proof of Theorem
\ref{maingen}. \hfill $\Box$

\section{Applications}
\label{applic}

\subsection{Intersecting sets with $\Bad(i_1,\ldots,i_N)$}
\label{badijo}

Let $\Bad(i_1,\ldots,i_N)$ be the set of $(i_1,\ldots,i_N)$-badly
approximable $N$-tuples in $\mathbb{R}^N$ as defined in
\S\ref{basicgen} and $\Bad(N) := \Bad(i_1,\ldots,i_N)$ with
$i_1=\ldots=i_N = 1/N$. Thus $\Bad(1)$ is simply the set $\Bad$ of
badly approximable real numbers.  Let $\Omega$ be a compact subset of
$\mathbb{R}^N$. The problem is to determine conditions on $\Omega$
under which
\begin{equation*}
  \Bad_{\Omega}(i_1, \ldots ,i_N) \ := \ \Omega \cap \Bad(i_1,\ldots,i_N)
\end{equation*}
is of full dimension; i.e. $ \dim \Bad_{\Omega}(i_1,\ldots,i_N) = \dim
\Omega $.  Recall, that the `$2$-dimensional' argument of
\S\ref{basicgen} can easily be extended to show that $ \dim
\Bad(i_1,\allowbreak \ldots,i_N) = N$.

To begin with, we address the above problem for the set
$\Bad_{\Omega}(N)  = \Omega \cap \Bad(N) $ in the case that
$\Omega$ supports an `absolutely $\alpha$-decaying' measure that
satisfies condition (A).

The notion of an `absolutely decaying' measure was introduced in
\cite{KLW}. The following restrictive definition, exploited in
\cite{PV1}, serves our purpose.  Let $\Omega$ be a compact subset of
$\mathbb{R}^N$ which supports a non-atomic, finite measure $m$.  Let
$\mathcal{L}$ denote a generic hyperplane of $\mathbb{R}^N$ and let
$\mathcal{L}^{(\epsilon)}$ denote its $\epsilon$-neighborhood. We say
that $m$ is \emph{absolutely $\alpha$-decaying} if there exist
strictly positive constants $C, \alpha, r_0$ such that for any
hyperplane $\mathcal{L}$, any $\epsilon > 0$, any $x \in \Omega$ and
any $r < r_0$,
\begin{equation*}
  m\left(B(x,r) \cap \mathcal{L}^{(\epsilon)} \right) \ \leq \ C \,
  \left(\dfrac{\epsilon}{r} \right)^{\alpha} m(B(x,r)) \ .
\end{equation*}
In the case $N=1$, the hyperplane $\mathcal{L}$ is simply a point $a
\in \mathbb{R}$ and $\mathcal{L}^{(\epsilon)}$ is the ball
$B(a,\epsilon)$ centred at $a$ of radius $\epsilon$. Also note that in
this case, if the measure $m$ satisfies condition (A) with exponent
$\delta$ then $m$ is automatically absolutely $\delta$-decaying.

\begin{thm}
  Let $\Omega$ be a compact subset of $\mathbb{R}^N$ which supports a
  measure $m$ satisfying condition (A) and which in addition is
  absolutely $\alpha$-decaying for some $\alpha > 0$.  Then
  \begin{equation*}
  \dim
  \Bad_{\Omega}(N) = \dim \Omega \ .
  \end{equation*}
  \label{thmab}
\end{thm}

\begin{pf*}{\textbf{Proof.}  \ }
  With reference to \S\ref{intro}, the set $\Bad_{\Omega}(N)$ can be
  expressed in the form $\Bad^*(\mathcal{R}, \beta, \rho)$ with
  $\rho(r) := r^{-(1+ \frac{1}{N})}$ and
  \begin{eqnarray*}
    & &X=(\mathbb{R}^N,d) \  , \ \  J:= \{ ((p_1,\ldots, p_N),q) \in
    \mathbb{N}^N \times \mathbb{N} \backslash \{0\}  \} \ , \ \
    \\ & & \\
    & & \alpha := ((p_1,\ldots, p_N),q) \in J \ , \ \
    \beta_{\alpha} := q \ , \ \ R_{\alpha}:= (p_1/q,\ldots , p_N/q) \ .
  \end{eqnarray*}
  Here $d$ is standard $\sup$ metric on $\mathbb{R}^N$; $d(x,y) :=
  \max \{ d(x_1,y_1), \ldots , d(x_N, y_N) \} $. Thus balls $B(c,r)$
  in $\mathbb{R}^N$ are genuinely cubes of sidelength $2r$.

  We show that the conditions of Theorem \ref{main} are satisfied.
  Clearly the function $\rho$ satisfies condition (B) and we are given
  that the measure $m$ supported on $\Omega $ satisfies condition (A).
  Also, since the resonant sets are points the condition that $\dim
  (\cup_{\alpha \in J} R_\alpha ) < \delta $ is satisfied. We need to
  establish the existence of the disjoint collection $\cthbn$ of balls
  (cubes) $2\theta B_{n+1}$ where $B_n$ is an arbitrary ball of radius
  $ k^{-n(1+ \tfrac{1}{N})}$ with centre in $\Omega$. In view of Lemma
  \ref{crucial}, there exists a disjoint collection $\cthbn$ such that
  \begin{equation}
    \# \cthbn \ \geq \ \kappa_1 \; k^{(1+ \frac{1}{N}) \delta } \ ;
    \label{absfrh1}
  \end{equation}
  i.e. \eqref{h1} of Theorem \ref{main} holds. We now verify that
  \eqref{h2} is satisfied for any such collection.

We consider two cases.

\noindent \emph{Case 1: $N=1$. \ } The
trivial argument of \S\ref{basiceg} shows that any interval $\theta
B_n $ with $\theta := \tfrac{1}{2} k^{-2}$ contains at most one
rational $p/q$ with $k^n \leq q < k^{n+1} $; i.e. $\alpha \in J(n+1)$.
Thus, for $k$ sufficiently large
\begin{equation*}
  \text{l.h.s. of \eqref{h2} } \leq \ 1 \ < \tfrac{1}{2} \times
  \text{r.h.s. of \eqref{absfrh1}} \ .
\end{equation*}
Hence (\ref{h2}) is trivially satisfied and Theorem \ref{main} implies
the desired result.

\noindent \emph{Case 2: $N \geq 2$. \ } We shall prove the theorem
in the case that $N= 2$.  There are no difficulties and no new ideas
are required in extending the proof to higher dimensions, especially
in view of the Simplex Lemma (see \S\ref{basicgen}).

Suppose that there are three or more rational points $(p_1/q,p_2/q) $
with $k^n \leq q < k^{n+1} $ lying within the ball/square $\theta B_n
$.  Now put $\theta := 2^{-1} (2 k^3)^{-1/2}$.  Then the `triangle'
argument of \S\ref{basicgen} (where $m$ is Lebesgue measure) implies
that the rational points must lie on a line $\mathcal{L}$ passing
through $\theta B_n$.  It follows that
\begin{alignat*}{2}
  \text{l.h.s. of \eqref{h2}} & \leq \# \left\{2\theta B_{n+1} \subset
    \cthbn : 2\theta B_{n+1} \cap \mathcal{L} \neq
    \emptyset \right\} \\[1em]
  & \leq \# \left\{2\theta B_{n+1} \subset \cthbn : 2\theta B_{n+1}
    \subset \mathcal{L}^{(\epsilon)} \right\} \quad \text{ for }
  \epsilon := 8 \theta
  k^{-(n+1)\tfrac{3}{2}} \\[1em]
  & \leq \dfrac{m(\theta B_n \cap \mathcal{L}^{(\epsilon)})}{m(
    2\theta B_{n+1}) } \hspace{33mm} \hfill \text{the \ balls $2\theta
    B_{n+1}$ are disjoint } \\[1em]
  &\leq a^{-1} b \, C \, \, 8^{\alpha} \, 2^{-\delta} \;
  k^{\frac{2}{3}(\delta - \alpha)} \hspace{21mm} m \text{ is
    absolutely $\alpha$-decaying} \\[1em]
  & < \tfrac{1}{2} \times \;\text{r.h.s. \ of \eqref{absfrh1}}
  \hspace{26.5mm} \hfill \text{for $k$ sufficiently large.}
\end{alignat*}
Hence (\ref{h2}) is satisfied  and Theorem \ref{main} implies the
desired result. \hfill $ \Box $
\end{pf*}

\vspace{3ex}

The following statement which combines Theorems 2.2 and 8.1 of
\cite{KLW}, shows that a large class of fractal measures are
absolutely $\alpha$-decaying and satisfy condition (A).

\vspace{3ex}

\begin{thm} \label{Theorem KLW}
  Let $\{\mathbf{S}_1, \dots, \mathbf{S}_k \} $ be an irreducible
  family of contracting self similarity maps of $\mathbb{R}^N$
  satisfying the open set condition and let $m$ be the restriction of
  $\mathcal{H}^{\delta}$ to its attractor $K$ where $ \delta := \dim
  K$.  Then $m$ is absolutely $\alpha$-decaying and satisfies
  condition (A).
\end{thm}

The simplest examples of such sets include regular Cantor sets, the
Sierpi\'nski gasket and the von Koch curve.  All the terminology
except for `irreducible' is pretty much standard -- see for example
\cite[Chp.9]{falc}. The notion of irreducible introduced in
\cite[\S2]{KLW} avoids the natural obstruction that there is a finite
collection of proper affine subspaces of $\mathbb{R}^N$ which is
invariant under $\{\mathbf{S}_1, \dots, \mathbf{S}_k \} $. More
recently, the class of examples regarding absolutely $\alpha$-decaying
measures has been extended by Urba{\'n}ski \cite{Urb1,Urb2}.

In view of Theorem \ref{Theorem KLW}, the following statement is a
simple consequence of Theorem \ref{thmab}. It has also been
independently established by Kleinbock \& Weiss \cite[Theorem
10.3]{KLW} and \cite{KW}. In fact, Theorem \ref{thmab} is also
derived in
\cite{KW} by an alternative approach.

\begin{cor}
  Let $\{\mathbf{S}_1, \dots, \mathbf{S}_k \} $ be an irreducible
  family of contracting self similarity maps of $\mathbb{R}^N$
  satisfying the open set condition and let $m$ be the restriction of
  $\mathcal{H}^{\delta}$ to its attractor $K$ where $ \delta := \dim
  K$.  Then
  \begin{equation*}
  \dim (K \cap \Bad(N)) = \dim K \ .
  \end{equation*}
  \label{corKW}
\end{cor}

We now consider the more general problem of determining conditions on
$\Omega$ under which $ \dim \Bad_{\Omega}(i_1,\ldots,i_N) = \dim
\Omega $. Under the hypotheses of Theorem \ref{maingen}, by modifying
the definition of `absolutely decaying' to accommodate `rectangles' it
is clearly possible to obtain an analogue of the `abstract' theorem
(Theorem \ref{thmab}) for $\Bad_{\Omega}(i_1, \ldots,i_N)$. We have
decided against establishing such a statement in this paper. The
reason for this is simple. We are currently unable to prove the
existence of a natural class of sets satisfying the more general
`rectangular' hypotheses. Nevertheless, in the special case that
$\Omega $ is a product space we are able to prove the following
statement.

\begin{thm}
  For $1\leq j \leq N$, let $\Omega_j$ be a compact subset of
  $\mathbb{R}$ which supports a measure $m_j$ satisfying condition (A)
  with exponent $\delta_j$. Let $\Omega$ denote the product set
  $\Omega_1 \times \ldots \times \Omega_N$.  Then, for any $N$-tuple
  $(i_1,\ldots, i_N)$ with $i_j \geq 0$ and $\sum_{j=1}^N \, i_j =1$,
  \begin{equation*}
    \dim \Bad_{\Omega}(i_1, ...,i_N) = \dim \Omega \ .
  \end{equation*}
\label{thm5}
\end{thm}

A simple application of the above theorem leads to following
result.

\begin{cor}
  Let $K_1$ and $K_2$ be regular Cantor subsets of $\mathbb{R}$. Then
  \begin{equation*}
    \dim \left( (K_1 \times K_2) \cap \Bad(i,j) \right) \ = \ \dim
    (K_1 \times K_2) \ = \ \dim K_1 \ + \ \dim K_2 \ .
  \end{equation*}
  \label{corthm5}
\end{cor}

\begin{pf*}{Proof of Theorem \ref{thm5}. \ }
  Without loss of generality assume that $N \geq 2$. The case that
  $N=1$ is covered by Theorem \ref{thmab}.  For the sake of clarity,
  as with the proof of Theorem \ref{thmab}, we shall restrict our
  attention to the case $N= 2$.

  Recall that since $\Omega_j \subset \mathbb{R}$ and $m_j$ satisfies
  (A), then $m_i$ is automatically absolutely $\delta_j$-decaying. A
  relatively straightforward argument shows that $m:=m_1\times m_2$ is
  absolutely $\alpha$-decaying on $\Omega$ with
  $\alpha:=\min\{\delta_1,\delta_2\}$. In fact this trivially follows
  from the following general fact - see \cite[\S9]{KLW}.

  \noindent \textbf{Fact:} {\em For $2\leq j \leq N$, if each $m_j$ is
  absolutely $\alpha_j$-decaying on $\Omega_j$, then
  $m:=m_1\times\ldots\times m_N$ is absolutely $\alpha$-decaying on
  $\Omega=\Omega_1\times\ldots\times\Omega_N$ with
  $\alpha=\min\{\alpha_1,\ldots,\alpha_N\}$. }

  Now let us write $\Bad(i,j)$ for $\Bad(i_1,i_2)$ and without loss of
  generality assume that $i<j$. The case $i=j$ is already covered by
  Theorem 4 since $m$ is absolutely $\alpha$-decaying on $\Omega$ and
  clearly satisfies condition (A). The set $\Bad_\Omega(i,j)$ can be
  expressed in the form $\Bad^*(\mathcal{R}, \beta, \rho_1,\rho_2)$
  with $\rho_1(r)=r^{-(1+i)}$, $\rho_2(r)=r^{-(1+j)}$ and
  \begin{eqnarray*}
    & &X=\mathbb{R}^2 \ , \ \  \Omega:= \Omega_1\times\Omega_2 \  , \
    \ J:= \{ ((p_1,p_2),q) \in \mathbb{N}^2 \times \mathbb{N}
    \setminus \{0\} \} \ , \ \
    \\[1em]
    & & \alpha := ((p_1,p_2),q) \in J \ , \ \ \beta_{\alpha} := q \ ,
    \ \ R_{\alpha}:= (p_1/q,p_2/q)   \ .
  \end{eqnarray*}
  With reference to Theorem \ref{maingencorr}, the functions
  $\rho_1,\rho_2$ satisfy condition (B*) and the measures $m_1, m_2$
  satisfy condition (A). Also note that $ \dim (\cup_{\alpha \in J}
  R_\alpha) = 0 $ since the union in question is countable.  We need
  to establish the existence of the collection $\cthfn$, where $F_n$
  is an arbitrary closed rectangle of size $2 k^{-n(1+i)} \times 2
  k^{-n(1+j)} $ with centre $c$ in $\Omega$. In view of Lemma
  \ref{crucial}, there exists a disjoint collection $\cthfn$ of
  rectangles $2\theta F_{n+1} \subset \theta F_n $ such that
  \begin{equation}
    \#\cthfn \ \geq \ \kappa_1 \, k^{(1+i)\delta_1}k^{(1+j)\delta_2} ;
    \label{cthfnlbd}
  \end{equation}
  i.e. \eqref{hgs1} of Theorem \ref{maingencorr} is satisfied. We now
  verify that \eqref{hgs2} is satisfied for any such collection.  With
  $\theta=2^{-1}(2k^3)^{-1/2}$, the `triangle' argument or
  equivalently the Simplex Lemma of \S\ref{basicgen} implies that
  \begin{equation}
    \text{l.h.s. of \eqref{hgs2}} \ \leq \ \#\{2\theta F_{n+1}
    \subset\cthfn: 2\theta F_{n+1} \cap \mathcal{L} \neq \emptyset\} \
    , \label{real}
  \end{equation}
  where $\mathcal{L}$ is a line passing through $\theta F_n$. Consider
  the thickening $T(\mathcal{L})$ of $\mathcal{L}$ obtained by
  placing rectangles $4\theta F_{n+1}$ centred at points of
  $\mathcal{L}$; that is, by `sliding' a rectangle $4\theta F_{n+1}$,
  centred at a point of $\mathcal{L}$, along $\mathcal{L}$.  Then,
  since the rectangles $2\theta F_{n+1}\subset\cthfn$ are disjoint,
  \begin{alignat}{2}
    \#\{2\theta F_{n+1}\subset\cthfn: 2\theta F_{n+1} &\cap\mathcal{L}
    \neq \emptyset\} \nonumber \\[1em]
    &\leq \#\{2\theta F_{n+1}\subset\cthfn: 2\theta F_{n+1} \subset
    T(\mathcal{L})\} \nonumber \\[1em]
    &\leq \ \dfrac{m(T(\mathcal{L})\cap\theta F_n)}{m(2\theta
      F_{n+1})}.
    \label{rhsij}
  \end{alignat}
  Without loss of generality we can assume that $\mathcal{L}$ passes
  through the centre of $\theta F_n$. To see this, suppose that
  $m(T(\mathcal{L})\cap\theta F_n)\neq 0$ since otherwise there is
  nothing to prove. Then, there exists a point $x \in
  T(\mathcal{L})\cap\theta F_n \cap \Omega$ such that
  \begin{equation*}
    T(\mathcal{L})\cap\theta F_n \subset 2\theta F_n^\prime \cap
    T^\prime(\mathcal{L}^\prime) \ .
  \end{equation*}
  Here $F_n^\prime $ is the rectangle of size $k^{-n(1+i)}\times
  k^{-n(1+j)}$ centred at $x$, $\mathcal{L}^\prime$ is the line
  parallel to $\mathcal{L}$ passing through $x$ and
  $T^\prime(\mathcal{L}^\prime)$ is the thickening obtained by
  `sliding' a rectangle $8\theta F_{n+1}$ centred at $x$, along
  $\mathcal{L}^\prime$. Then the following argument works just as well
  on $2\theta F_n^\prime \cap T^\prime(\mathcal{L}^\prime)$.

  Let $\Delta$ denote the slope of the line $\mathcal{L}$ and assume
  that $\Delta \geq 0$. The case $\Delta < 0$ can be dealt with
  similarly.  By moving the rectangle $\theta F_n$ to the origin,
  straightforward geometric considerations lead to the following
  facts:

\begin{enumerate}[(F1)]
   \item
  \begin{equation*}
  T(\mathcal{L})= \mathcal{L}^{(\epsilon)} \quad \text{where}
      \quad \epsilon:= \dfrac{4\theta\left(k^{-(n+1)(1+j)}+\Delta
          k^{-(n+1)(1+i)}\right)}{\sqrt{1+\Delta^2}}
    \end{equation*}
  \item $T(\mathcal{L})\cap\theta F_n \subset F(c;l_1, l_2)$ where
    $F(c;l_1, l_2)$ is the rectangle with the same centre $c$ as $F_n$
    and of size $2l_1\times 2l_2$ with
    \begin{eqnarray*}
      l_1& := & \  \dfrac{\theta}{\Delta} \left(k^{-n(1+j)} + 4k^{-(n+1)(1+j)}
        + \Delta k^{-(n+1)(1+i)}\right)  \\[1em]
      l_2 & := & \ \theta k^{-n(1+j)} \  \ .
    \end{eqnarray*}
 \end{enumerate}
  We now estimate the right hand side of (\ref{rhsij}) by considering
  two cases.  Throughout, let $a_i, b_i$ denote the constants
  associated with the measure $m_i$ and condition (A) and let
  \begin{equation*}
  \varpi \, := \, 3 \; \left(\dfrac{4b_1b_2}{\kappa_1 \,
      a_1a_22^{\delta_1+\delta_2}}\right)^{1/\delta_1} \ .
  \end{equation*}

  \noindent \textbf{Case (i): }  $\Delta \geq \varpi
  k^{-n(1+j)}/k^{-n(1+i)}$. \ In view of (F2) above, we trivially have
  that
  \begin{equation*}
    m(\theta F_n \cap T(\mathcal{L})) \ \leq  \ m(F(c;l_1, l_2)) \
    \leq \ b_1\, b_2 \, l_1^{\delta_1} \, l_2^{\delta_2} \ .
  \end{equation*}
  It follows that
  \begin{alignat*}{2}
    \dfrac{m(T(\mathcal{L})\cap\theta F_n)}{m(2\theta F_{n+1})} &\leq
    \dfrac{b_1b_2 l_1^{\delta_1} l_2^{\delta_2}}{a_1a_2
      (2\theta)^{\delta_1+\delta_2} \, k^{-(n+1)(1+j)\delta_1} \,
      k^{-(n+1)(1+i)\delta_2}} \\[1em]
    & \leq \dfrac{b_1b_2}{a_1a_22^{\delta_1+\delta_2}}
    \left(\dfrac{1}{\varpi} +\dfrac{1}{\varpi k^{1+j}}+
      \dfrac{1}{k^{1+i}} \right)^{\delta_1} \ k^{(1+j)\delta_1} \,
    k^{(1+i)\delta_2} \\[1em]
    & \leq \dfrac{b_1b_2}{a_1a_22^{\delta_1+\delta_2}}
    \left(\dfrac{3}{\varpi} \right)^{\delta_1} \ k^{(1+j)\delta_1} \,
    k^{(1+i)\delta_2}  \\[1em] & =  \dfrac{\kappa_1}{4} \ k^{(1+j)\delta_1}
    k^{(1+i)\delta_2} \ .
  \end{alignat*}

  \noindent \textbf{Case (ii): }  $0 \leq\Delta < \varpi
  k^{-n(1+j)}/k^{-n(1+i)}$. By the covering lemma of \S\ref{prel},
  there exists a collection $\mathcal{B}_n$ of disjoint balls $B_n$
  with centres in $\theta F_n\cap\Omega$ and radii $\theta
  k^{-n(1+j)}$ such that
  \begin{equation*}
    \theta F_n \cap \Omega \ \subset \  \bigcup_{B_n \in
      \mathcal{B}_n} 3B_n \ .
  \end{equation*}
  Since $i<j$, it is easily verified that the disjoint collection
  $\mathcal{B}_n$ is contained in $2\theta F_n$ and thus
  $\#\mathcal{B}_n \leq m(2\theta F_n)/m(B_n)$. It follows that
  \begin{alignat*}{3}
    m(\theta F_n &\cap \, T(\mathcal{L})) \leq m\left( \cup_{B_n \in
        \mathcal{B}_n} 3B_n \cap \; T(\mathcal{L})\right) & \\[1em]
    &\leq \#\mathcal{B}_n \  m(3B_n \cap T(\mathcal{L})) & \\[1em]
    &\leq \dfrac{m(2\theta F_n)}{m(B_n)} \ m\left(3B_n \cap
      \mathcal{L}^{(\epsilon)}\right) & \text{by (F1) above} \\
    &\leq m(2\theta F_n) \, \dfrac{m(3B_n)}{m(B_n)} \
    \left(\dfrac{\epsilon}{3\theta k^{-n(i+j)}}\right)^\alpha & \quad m
    \text{ is absolutely $\alpha$-decaying.}
  \end{alignat*}

  Now notice that
  \begin{equation*}
    \dfrac{\epsilon}{3\theta k^{-n(i+j)}} \ \leq \  \dfrac{4}{3} \,
    (k^{-(1+j)} + \varpi k^{-(1+i)}) \ .
  \end{equation*}
  Hence, for $k$ sufficiently large we have that
  \begin{equation*}
    \dfrac{m(T(\mathcal{L})\cap\theta F_n)}{m(2\theta F_{n+1})} \ \leq
    \ \dfrac{\kappa_1}{4} \ k^{(1+j)\delta_1} k^{(1+i)\delta_2} \ .
  \end{equation*}

  On combining the above two cases, we have that
  \begin{equation*}
    \text{l.h.s. of \eqref{hgs2}} \leq \dfrac{m(T(\mathcal{L}) \cap
      \theta F_n)}{m(2\theta F_{n+1})} \ \leq \ \dfrac{\kappa_1}{4}  \
    k^{(1+j)\delta_1} k^{(1+i)\delta_2} \ = \ \dfrac{1}{4} \times  \
    \text{l.h.s. of \eqref{cthfnlbd}} \ .
  \end{equation*}
  Hence \eqref{hgs2} is satisfied and Theorem \ref{maingencorr}
  implies the desired result.

\hfill $ \Box $

\end{pf*}

The argument used to establish Theorem \ref{thm5} can be adapted
in the obvious manner to prove a slightly more general result.

\begin{thm}
  For $1\leq j \leq N$, let $\Omega_j $ be a compact subset of
  $\mathbb{R}^{d_j}$ which supports an absolutely $\alpha_j$-decaying
  measure $m_j$ satisfying condition (A) with exponent $\delta_j$. Let
  $\Omega$ denote the product set $\Omega_1 \times \ldots \times
  \Omega_N$.  Then, for any $N$-tuple $(i_1,\ldots, i_N)$ with $i_j
  \geq 0$ and $\sum_{j=1}^N d_j \, i_j =1$,
  \begin{equation*}
    \dim \Bad_{\Omega}( \, \underbrace{i_1, \ldots, i_1}_{\text{$d_1$
        times}}\, ;\ \underbrace{i_2, \ldots, i_2}_{\text{$d_2$
        times}} \, ; \ \ldots ;\ \underbrace{i_N\, , \ldots,
      i_N}_{\text{$d_N$ times}}\;) \ = \ \dim \Omega \ = \
    \sum_{j=1}^N \delta_j \ .
  \end{equation*}
  \label{thm6}
\end{thm}

The following is a simple consequence of Theorem \ref{Theorem KLW}
and Theorem \ref{thm6}.

\begin{cor}
  For $1\leq j \leq N$, let $K_j $ be the attractor of a finite
  irreducible family of contracting self similarity maps of
  $\mathbb{R}^{d_j}$ satisfying the open set condition.  Let $m_j$ be
  the restriction of $\mathcal{H}^{\delta_j}$ to $K_j$ where
  $\delta_j= \dim K_j$. Let $K$ denote the `product attractor' $K_1
  \times \ldots \times K_N$.  Then, for any $N$-tuple $(i_1,\ldots,
  i_N)$ with $i_j \geq 0$ and $\sum_{j=1}^N d_j \, i_j =1$,
  \begin{equation*}
  \dim (K \, \cap \, \Bad(\, \, \underbrace{i_1, \ldots,
    i_1}_{\text{$d_1$ times}}\, ;\ \underbrace{i_2, \ldots,
    i_2}_{\text{$d_2$ times}}\, ; \ \ldots ;\ \underbrace{i_N\, ,
    \ldots, i_N}_{\text{$d_N$ times}}\;)) \ = \ \dim K  \ .
  \end{equation*}
  \label{attractprod}
\end{cor}

As an application of Corollary \ref{attractprod} we obtain the
following statement which to some extent is more illuminating -- even
this special case appears to be new.

\begin{cor}
  Let $V \subset \mathbb{R}^2$ be the von Koch curve and $K \subset
  \mathbb{R}$ be the middle third Cantor set. Then, for any positive
  $i$ and $j$ with $2 \, i + j =1$
  \begin{equation*}
    \dim \left( (V \times K) \cap \Bad(i,i,j) \right) \ = \ \dim (V
    \times K)  \ = \ \dfrac{\log 8}{\log 3}  \ .
  \end{equation*}
  \label{corthm6}
\end{cor}

\subsubsection{Remarks related to Schmidt's conjecture}
\label{schsec}

In \S\ref{basicgen}, we mentioned the result that $\dim
(\Bad(i,j)\cap\Bad(1,0)\cap\Bad(0,1)) = 2$. This can easily be
obtained via Theorem \ref{thm5}. To see this, first of all notice that
$\Bad \times \Bad = \Bad(1,0) \cap \Bad(0,1)$. For $M \geq 2$, let
$F_M :=\{ x \in [0,1]: x:=[a_1,a_2, \ldots] \text{ with } a_i \leq M
\text{ for all } i\}$. Thus $F_M$ is the set of real numbers in the
unit interval with partial quotients bounded above by $M$. By
definition $F_M$ is a compact subset of $\Bad$ and moreover it is well
known that $F_M$ supports a measure $m_M$ which satisfies condition
(A) with exponent $\delta_M$ with $\delta_M \rightarrow 1$ as $M
\rightarrow \infty$. Now let $\Omega:= F_M \times F_M$, then Theorem
\ref{thm5} implies that
\begin{equation*}
  \dim(\Bad(i,j) \cap \Bad(1,0) \cap \Bad(0,1)) \ \geq \
  \dim(\Bad_\Omega(i,j)) \ = \ 2\delta_M  \ .
\end{equation*}
On letting $M \rightarrow \infty$, we obtain that $\dim(\Bad(i,j) \cap
\Bad(1,0) \cap \Bad(0,1)) \geq 2$. The complementary upper bound
result is trivial since the set in question is a subset of
$\mathbb{R}^2$.

Recall, that Schmidt's conjecture \cite{Schmidt} states that
$\Bad(i,j) \cap \Bad(i',j') \neq \emptyset $. In fact, Schmidt stated
this conjecture in the simpler situation when $i=j'=1/3$ and
$i'=j=2/3$. Even this specific and symmetric case is unsolved.  In
order to illustrate a possible approach towards the conjecture via the
results of this paper we consider the special case of $ \Bad(i,j) \cap
\Bad(1/2,1/2) $.  Suppose for the moment that we could find a compact
set $\Omega \subseteq \Bad(i,j)$ with a measure $m$ satisfying
condition (A) for some $\delta > 1$. Let $\rho(r) = r^{-3/2}$. Using
Lemma \ref{crucial} together with the `triangle' argument or
equivalently the Simplex Lemma of \S\ref{basicgen}, we may construct
collections $\cthbn$ as in the statement of Theorem \ref{main}. The
condition that $\delta > 1$ is used to ensure \eqref{h2}. This leads
to the following enticing statement:

\emph{If there exists a compact subset $\Omega$ of $ \Bad(i,j) $ which
  supports a measure $m$ satisfying condition (A) with exponent
  $\delta > 1 $, then $$ \dim ( \Bad(i,j) \cap
  \Bad(1/2,1/2)) \geq \delta  \ . $$}

Clearly, this would imply that $ \Bad(i,j) \cap \Bad(1/2,1/2) \neq
\emptyset $. Regarding the above statement, it is not particularly
difficult to prove the existence of a compact subset $\Omega$
supporting a measure $m$ satisfying condition (A) with $\delta < 1 $.
However, from this we are not able to deduce that $ \dim ( \Bad(i,j)
\cap \Bad(1/2,1/2)) \geq \delta $ or even that $ \Bad(i,j) \cap
\Bad(1/2,1/2) \neq \emptyset $.

\subsection{Rational Maps}

In this section we consider the `badly approximable' analogue of the
`shrinking target' problem introduced in \cite{inv} for expanding
rational maps.  Let $T$ be an expanding rational map (degree $\geq 2$)
of the Riemann sphere $\overline{\mathbb{C}} = \mathbb{C} \cup \{
\infty \}$ and $J(T)$ be its Julia set. For any $z_0 \in J(T)$
consider the set
\begin{multline*}
  \Bad_{z_0}(J) \, := \, \{ z \in J(T)\ : \exists \; c(z)>0 \text{
    such that } \\
  T^n(z) \notin B \, ( z_0 , c(z) ) \text{ for any } n \in
  \mathbb{N}\} \ .
\end{multline*}
Clearly, the forward orbit of points in $\Bad_{z_0}(J)$ are not dense
in $J(T)$.  Now let $m$ be Sullivan measure and $\delta = \dim J(T)$.
Thus $m$ is a non-atomic, $\delta$-conformal probability measure
supported on $J(T)$ and since $T$ is expanding it satisfies condition
(A). Moreover, $m$ is equivalent to $\delta$-dimensional Hausdorff
measure $\mathcal{H}^{\delta}$ -- see \cite{inv,ihes} for the
details.  In view of the `Khintchine type' result for expanding
rational maps (see, for example \cite[\S8.4]{BDV}) it is easily
verified that $\mathcal{H}^{\delta}(\Bad_{z_0}(J)) = 0 =
m(\Bad_{z_0}(J)) $. Nevertheless, the set $\Bad_{z_0}(J)$ is large in
that it is of maximal dimension.

\begin{thm}
  \begin{equation*}
    \dim \Bad_{z_0}(J)  \ =  \ \delta \ .
  \end{equation*}
 \label{thmrat}
\end{thm}
This result is not new and has been established by numerous people
(see e.g. \cite{MR1444052}).  However, we give a short proof which
indicates the versatility and generality of our framework and results.

\begin{pf*}{Proof of Theorem \ref{thmrat}. \ }
  In view of the bounded distortion property for expanding maps
  \cite[Proposition 1]{inv}, we can rewrite $\Bad_{z_0}(J)$ in terms
  of points in the Julia set which `stay clear' of balls centred
  around the backward orbit of the selected point $z_0$:
  \begin{multline*}
    \Bad_{z_0}(J)  \equiv  \{ z \in J(T)\, :  \exists \; c(z)>0 \text{
      such that } \\ z \notin  B \, (y, c(z)|(T^n)^\prime(y)|^{-1}) \text{ for
      any } (y,n) \in I \} \ ,
  \end{multline*}
  where $I := \{(y,n) : n \in \mathbb{N} \text{ with } T^n(y) = z_0 \}
  $. Also, since $T$ is expanding, $J(T)$ can be thought of as a
  compact metric space with the usual metric on $\mathbb{C}$. It is
  now clear that $\Bad_{z_0}(J)$ can be expressed in the form
  $\Bad^*(\mathcal{R}, \beta, \rho)$ with $\rho(r):=r^{-1}$ and
  \begin{equation*}
    X = \Omega := J(T)  \  , \ \ J:= I \ , \ \ \alpha := (y,n) \in
    I \ , \ \ \beta_\alpha := |(T^n)^\prime(y)| \ , \ \ R_{\alpha}:= y
    \ .
  \end{equation*}
  With reference to Theorem \ref{main}, Sullivan measure $m$ and the
  function $\rho$ satisfy condition (A) and (B) respectively. To
  deduce Theorem \ref{thmrat} from Theorem \ref{main} we need to
  establish the existence of the disjoint collection $\cthbn$ of balls
  $2\theta B_{n+1}$ where $B_n$ is an arbitrary ball of radius $
  k^{-n}$ with centre in $\Omega$. In view of Lemma \ref{crucial}, for
  $k$ sufficiently large, there exists a disjoint collection $\cthbn$
  such that
  \begin{equation}
    \# \cthbn \ \geq \ \kappa_1 \; k^{ \delta } \ ; \label{absfrh2}
  \end{equation}
  i.e. \eqref{h1} of Theorem \ref{main} holds.  We now verify that
  \eqref{h2} is satisfied for any such collection. First we recall a
  key result which is the second part of the statement of Lemma 8 in
  \cite{ihes}. For ease of reference we keep the same notation and
  numbering of constants as in \cite{ihes}.

  \noindent \emph{Constant Multiplicity:} For $X \in \mathbb{R}^+$, let
  $P(X)$ denote the set of pairs $(y,n) \in I$ such that $ f_n(y) -
  C_8 \leq X \leq f_{n+1}(y) + C_8 $, where $f_n(y) := \log
  |(T^n)^\prime(y)|$.  Let $z \in J(T)$. Then there are no more than
  $C_9$ pairs $(y,n) \in P(X)$ such that $z \, \in \, B \left( y,
    C_{10} \; |(T^n)^\prime(y)|^{-1} \right) $.

  We are now in the position to verify (\ref{h2}) of Theorem
  \ref{main}.  By definition $J(n+1):=\{(y,m)\in I: k^{n-1} \leq \vert
  (T^m)^\prime(y) \vert < k^{n} \}$ and let $\theta:= C_{10}k^{-1}$.
  It follows that
  \begin{alignat}{2}
    \text{l.h.s. of \eqref{h2}} & \leq \#\{y \in \theta B_n :
    (y,m)\in J(n+1)\}  \nonumber  \\[1em]
    &\leq \#\{y\in B(c, C_{10}\vert (T^m)^\prime(y)\vert^{-1}):
    (y,m)\in J(n+1)\} \ ,
    \label{shit}
  \end{alignat}
  where $c$ is the centre of $\theta B_n $. Without loss of
  generality, assume that $ |T^\prime(z_0)|>1$.  Otherwise, since $T$
  is expanding we simply work with some higher iterate $T^q$ of $T$
  for which $ |(T^q)^\prime (z_0)| > 1 $. Then, the chain rule
  together with the above `constant multiplicity' fact implies that
  the r.h.s. of \eqref{shit} is $\ll C_9 \log k $. Hence, for $k$
  sufficiently large
  \begin{equation*}
    \text{l.h.s. of \eqref{h2}} \ \leq \ \tfrac{1}{2} \times \;
    \text{r.h.s. of \eqref{absfrh2}} \ .
  \end{equation*}
  Thus, \eqref{h2} is easily satisfied and Theorem \ref{main} implies
  Theorem \ref{thmrat}. \hfill $\Box$
\end{pf*}

\noindent \textbf{Remark: \; }  It is worth mentioning that  our
framework also yields (just as easily) the analogue of Theorem
\ref{thmrat} within the Kleinian group setup.  Briefly, let $G$ be
either a geometrically finite Kleinian group of the first kind or a
convex co-compact group and let $\Lambda(G)$ denote its limit set. For
these groups, Patterson measure supported on $\Lambda(G)$ satisfies
condition (A) and plays the role of Sullivan measure.  Then, it is not
difficult to obtain the Kleinian group analogue of Theorem
\ref{thmrat} via Theorem \ref{main}; i.e. the set of `badly
approximable' limit points is of full dimension -- $\dim \Lambda(G)$.

\subsection{Complex numbers}

In this section we consider the badly approximable analogue of
$\Bad(i_1,\ldots,i_N)$ in $\mathbb{C}^N$. Let $N \in \mathbb{N}$ and
$i_1, \ldots, i_N \geq 0$ such that $i_1 + \cdots + i_N = 1$.  Now
define the set $\Bad_{\mathbb{C}}(i_1,\ldots, i_N)$ to consist of $z
:=(z_1, \ldots , z_N) \in \mathbb{C}^N $ for which there exists a
constant $c(z) > 0$ such that for any $q, p_1, \dots, p_N \in
\mathbb{Z}[i]$, $q \neq 0$,
\begin{equation*}
  \max\{ \vert q z_1 - p_1 \vert^{1/i_1}, \dots, \vert q z_N - p_N
  \vert^{1/i_N} \} \geq c(z) \vert q \vert^{-1} \ .
\end{equation*}
In the case $i_1= \cdots = i_N = 1/N$, the corresponding set will
be denoted by $\Bad_{\mathbb{C}}(N)$.  Notice, that the role of the
rationals in the real setup is replaced by ratios of Gaussian integers
in the complex setup. We shall refer to the latter as Gaussian points.

The Hausdorff dimension of the set $\Bad_{\mathbb{C}}(N)$ has been
studied in the past by various people using Kleinian groups
\cite{MR98k:22043}, Riemannian geometry \cite{MR95m:30059} and
Schmidt's $(\alpha, \beta)$-games \cite{dodson:_hausd_dioph}.  Theorem
\ref{main} of this paper will also give the Hausdorff dimension of
this set. In fact, our general framework enables us to find the
dimension of $\Bad_{\mathbb{C}}(i_1,\ldots, i_N)$ intersected with
direct products of sets supporting measures satisfying condition (A).
As a consequence, the previously known results are extended to the
`rectangular' or `weighted' form of simultaneous approximation in
$\mathbb{C}^N$. The following statement is the `complex' analogue of
Theorem \ref{thm5}.

\begin{thm}
  For $1\leq j \leq N$, let $\Omega_j $ be a compact subset of
  $\mathbb{C}$ which supports a measure $m_j$ satisfying condition (A)
  with exponent $\delta_j$. Let $\Omega$ denote the product set
  $\Omega_1 \times \ldots \times \Omega_N$.  Then, for any $N$-tuple
  $(i_1,\ldots, i_N)$ with $i_j \geq 0$ and $\sum_{j=1}^N \, i_j =1$,
  \begin{equation*}
  \dim (\Bad_{\mathbb{C}}(i_1,\ldots, i_N) \cap \Omega \, ) =
  \dim \Omega \ .
  \end{equation*}
  \label{thmcomplex}
\end{thm}

The following complex notion of absolutely decaying measures will be
useful in proving the above theorem. Let $\Omega$ be a compact subset
of $\mathbb{C}^N$ which supports a non-atomic, finite measure $m$. Let
$\mathcal{L}$ denote a generic $(N-1)$-dimensional complex hyperplane
of $\mathbb{C}^N$ and let $\mathcal{L}^{(\epsilon)}$ denote its
$\epsilon$-neighborhood. We say that $m$ is \emph{absolutely
  $\alpha$-decaying} if there exist strictly positive constants $ C,
\alpha, r_0 $ such that for any complex hyperplane $\mathcal{L}$, any
$\epsilon > 0$, any $z \in \Omega$ and any $r < r_0$,
\begin{equation*}
  m\left(B(z,r) \cap \mathcal{L}^{(\epsilon)} \right) \ \leq \ C \,
  \left(\dfrac{\epsilon}{r} \right)^{\alpha} m(B(z,r)) \ .
\end{equation*}
Note that if $N=1$, so that $\Omega $ is a subset of $\mathbb{C}$, the
complex hyperplane $\mathcal{L}$ is simply a point $a \in \mathbb{C}$
and $\mathcal{L}^{(\epsilon)}$ is the ball $B(a,\epsilon)$ centred at
$a$ of radius $\epsilon$.  Moreover, if the measure $m$ satisfies
condition (A) with exponent $\delta$ then $m$ is automatically
absolutely $\delta$-decaying.

It is easy to verify that the statement of the `Fact' in
\S\ref{badijo} regarding the product of absolutely decaying measures
remains valid for the complex notion.

\begin{pf*}{Proof of Theorem \ref{thmcomplex} \emph{(Sketch}). \ }
  As usual we restrict our attention to the case $N = 2$ and write
  $\Bad_{\mathbb{C}}(i,j)$ for $\Bad_{\mathbb{C}}(i_1,i_2)$. Assume
  that $i \leq j$. Clearly, the set $\Bad_{\mathbb{C}}(i,j) \cap
  \Omega $ can be expressed in the form $\Bad^*(\mathcal{R}, \beta,
  \rho_1,\rho_2)$ with $\rho_1(r)=r^{-(1+i)}$, $\rho_2(r)=r^{-(1+j)}$
  and
  \begin{eqnarray*}
    & &X= (\mathbb{C}^2,d)  \
    , \ \  J:= \{ ((p_1, p_2) ,q) \in \mathbb{Z}[i]^2 \times
    \mathbb{Z}[i] \backslash \{0\}\} \ , \\
    & & \\ & & \alpha := ((p_1, p_2),q) \in J \ , \ \ \beta_{\alpha} := \vert
    q \vert \ , \ \ R_{\alpha}:=  (p_1/q, p_2,q) \ .
  \end{eqnarray*}
  The metric $d$ on $\mathbb{C}^2$ is the maximum of the coordinate
  metrics; i.e. $d((z_1, z_2), \allowbreak (z_1', z_2')) = \max
  \{d(z_1, z_1'), d(z_2, z_2')\}$.  Also note that the measure
  $m:=m_1\times m_2$ is absolutely $\alpha$-decaying on $\Omega$ with
  $\alpha:=\min\{\delta_1,\delta_2\}$. This follows from the above
  discussion concerning the complex notion of absolutely decaying
  measures and their product.

  With reference to Theorem \ref{maingencorr}, we need to establish
  the existence of the collection $\cthfn$ where $F_n $ is an
  arbitrary closed polydisc $B_{n,1} \times B_{n,2}$ with centre $c$
  in $\Omega$. Here $B_{n,1}$ (resp. $B_{n,2}$) is a closed ball in
  $\mathbb{C}$ of radius $ k^{-n(1+i)} $ (resp. $k^{-n(1+j)} $). In
  view of Lemma \ref{crucial}, there exists a disjoint collection
  $\cthfn$ of polydiscs $2\theta F_{n+1} \subset \theta F_n $ such
  that \eqref{hgs1} of Theorem \ref{maingencorr} is satisfied. We now
  verify that (\ref{hgs2}) is satisfied for any such collection by
  modifying the proof of Theorem \ref{thm5} in the obvious manner.
  The only part which is not so obvious is the complex analogue of the
  `triangle' argument of \S\ref{basicgen}. For this suppose that
  $\theta F_n$ is given and that there are at least three Gaussian
  points $(p_1/q,p_2/q), (p'_1/q',p'_2/q') $ and
  $(p''_1/q'',p''_2/q'') $ with
  \begin{equation*}
    k^n \leq |q|,|q'|,|q''| < k^{n+1}
  \end{equation*}
  lying within $\theta F_n $. Suppose for the moment that they do not
  lie on a one--dimensional complex hyperplane (i.e. a complex line)
  $\mathcal{L}$ of $\mathbb{C}^2$ and consider the determinant
  \begin{displaymath}
    D = \det
    \begin{pmatrix}
      1 & p_1/q & p_2/q \\
      1 & p_1'/q' & p_2'/q' \\
      1 & p_1''/q'' & p_2''/q'' \\
    \end{pmatrix} \neq 0 \ \ .
  \end{displaymath}

  Expanding the determinant in the first column and using the fact
  that the ring of Gaussian integers is a unique factorization domain,
  we find that
  \begin{displaymath}
    \vert D \vert > \dfrac{1}{k^{3(n+1)}} \ \ .
  \end{displaymath}
  On the other hand, the absolute value of $D$ can be at most twice
  the diameters of the two projections $\theta B_{n,1}$ and $\theta
  B_{n,2}$ of $\theta F_n$. That is
  \begin{displaymath}
    \vert D \vert \ \leq \  2  \; \dfrac{2\, \theta}{k^{n(i+1)}}
    \, \dfrac{2 \, \theta}{k^{n(j+1)}} \ = \ \dfrac{8
      \theta^2}{k^{3n}} \ \ .
  \end{displaymath}

  To see this, note that for $(z_1, z_2), (z_1', z_2'), (z_1'', z_2'')
  \in \theta F_n$
  \begin{eqnarray*}
    \left \vert \det
      \begin{pmatrix}
        1 & & z_1 &  & z_2 \\
        1 & & z_1' & & z_2' \\
        1 & & z_1'' & & z_2'' \\
      \end{pmatrix}
    \right \vert & \ = \  &  \vert (z_1 - z_1')(z_2' - z_2'') + (z_1' -
    z_1'')(z_2' - z_2) \vert \\ & & \\ & \ \leq \ & 2 \times 2 \theta \rho_1(k^n) \,  2  \theta
    \rho_2(k^n) \ .
  \end{eqnarray*}
  Now with $\theta := (8 k^3)^{-1/2} $, we obtain the desired
  contradiction. Thus, if there are two or more Gaussian points with
  $k^n \leq |q| < k^{n+1} $ lying within $\theta F_n$ then they must
  lie on a complex line $\mathcal{L}$. It now follows that
  \begin{equation*}
    \text{l.h.s. of \eqref{hgs2}} \ \leq \ \#\{2\theta F_{n+1}
    \subset\cthfn: 2\theta F_{n+1} \cap \mathcal{L} \neq \emptyset\} \
    .
  \end{equation*}
  This is the precise complex analogue of \eqref{real} and the proof
  can now be completed by modifying the proof of the real case
  (Theorem \ref{thm5}) in the obvious manner. We leave the details to
  the reader. \hfill $ \Box $
\end{pf*}

\medskip

It is worth mentioning that Theorem \ref{thmcomplex} can be
generalized in the obvious manner to obtain the complex analogue
of Theorem \ref{thm6}.

\subsection{$p$-adic numbers}

For a prime $p$, let $| \cdot |_p$ denote the $p$-adic absolute value
and let $\mathbb{Q}_p$ denote the $p$-adic field.  Furthermore, let
$\mathbb{Z}_p$ denote the ring of $p$-adic integers. In this section
we consider the badly approximable analogue of $\Bad(i_1,\ldots,i_N)$
in $\mathbb{Z}_p^N$. Let $N \in \mathbb{N}$ and $i_1, \ldots, i_N \geq
0$ such that $i_1 + \cdots + i_N = 1$.  Now define the set
$\Bad_{\mathbb{Z}_p}(i_1,\ldots, i_N)$ to consist of $x :=(x_1,
\ldots, x_N) \in \mathbb{Z}_p^N$ for which there exists a constant
$c(x) > 0$ such that
\begin{eqnarray}
  \max\{\vert q x_1 - r_1 \vert_p^{1/(1+i_1)}, \ldots & , &  \vert q
  x_N -  r_N \vert_p^{1/(1+i_N)} \} \nonumber \\[0.3em]  & \geq &  \ c(x) \, \max\{ \vert r_1 \vert
  , \ldots, \vert r_N \vert, \vert q \vert \}^{-1} \ ,
  \label{psetup}
\end{eqnarray}
for all $ ((r_1, \ldots, r_N) ,q) \in \mathbb{Z}^N \times
\mathbb{Z}\setminus \{0\} $.  In the case $i_1= \cdots = i_N = 1/N$,
the corresponding set will be denoted by $\Bad_{\mathbb{Z}_p}(N)$.

There are two points worth making when comparing the above set with
the `classical' set $\Bad(i_1, ...,i_N)$. Firstly, the r.h.s of
\eqref{psetup} in the $p$-adic setup is a function of $\max(|r_1|,
\ldots, |r_N|,|q|)$ rather than simply $|q|$.  This is due to the fact
that within the $p$-adic setup for any $x \in \mathbb{Z}_p^N$ and
$q\in\mathbb{Z}$ there exists $r \in \mathbb{Z}^N$ such that the l.h.s
of \eqref{psetup} can be made arbitrarily small. Thus, the set of $x
\in \mathbb{Z}_p^N$ for which l.h.s of (\ref{psetup}) $ \geq c(x) \,
\vert q \vert^{-1} $ is in fact empty and there is nothing to prove.
Secondly, in the $p$-adic setup the `weighting' factor occurring on
the l.h.s of (\ref{psetup}) is $1/(1+i_s)$ rather than $1/i_s$ ($1\leq
s \leq N$). This is due to the fact that we approximate in terms of
the $p$-adic absolute value on the left hand side, but measure the
`rate' of approximation in terms of the ordinary absolute value on the
right hand side. Because of the arithmetical properties of the p-adic
absolute value, we generally expect the `rate' of the approximation to
be better (see below).

Badly approximable $p$-adic numbers have in the past been studied by
Abercrombie \cite{MR96g:11078}, who showed that
$\Bad_{\mathbb{Z}_p}(1)$ has full Hausdorff dimension. In higher
dimensions, the corresponding result for even the `symmetric' set
$\Bad_{\mathbb{Z}_p}(N)$ is unknown.  Using the framework established
in this paper, we are able to prove the following complete result.

\begin{thm}
  \begin{equation*}
    \dim \Bad_{\mathbb{Z}_p}(i_1,\dots, i_N) = N \ .
  \end{equation*}
  \label{thmpadic}
\end{thm}

\begin{pf*}{Proof of Theorem \ref{thmpadic} \emph{(Sketch)}. \ }
  As in the preceding applications, we restrict our attention to the
  case $N = 2$ and write $\Bad_{\mathbb{Z}_p}(i,j)$ for
  $\Bad_{\mathbb{Z}_p}(i_1,i_2)$.  Assume that $i \leq j$. Clearly the set
  $\Bad_{\mathbb{Z}_p}(i,j)$ can be expressed in the form
  $\Bad^*(\mathcal{R}, \allowbreak \beta, \rho_1,\rho_2)$ with $\rho_1(x) :=
  x^{-(1+i)}$, $\rho_2(x) := x^{-(1+j)}$ and
  \begin{eqnarray*}
    & &X= \Omega:= \mathbb{Z}_p^2 = \mathbb{Z}_p \times \mathbb{Z}_p \
    , \ \  J:= \{ ((r_1, r_2) ,q) \in \mathbb{Z}^2 \times \mathbb{Z}
    \setminus \{0\} \} \ , \\[0.6em]
    & & \alpha := ((r_1, r_2),q) \in J \ , \ \ \beta_{\alpha} := \max\{
    \vert r_1 \vert, \vert r_2 \vert , \vert q \vert \} \ , \\[0.6em]
    & & R_{\alpha}:=  \{(x_1,x_2) \in \mathbb{Z}_p^2 : qx_1 = r_1, \ q
    x_2 = r_2 \} \ \ .
  \end{eqnarray*}
  Furthermore, $d = d_1 \times d_1 $ where $d_1(x,y):=|x-y|_p$ is the
  $p$-adic metric on $\mathbb{Q}_p$ and $ m:= \mu \times \mu$ where
  $\mu$ is normalized Haar measure on $\mathbb{Q}_p$. Thus,
  $\mu(\mathbb{Z}_p)=1$ and $\mu( B(x,p^{-t})) = p^{-t}$ for any $t
  \in \mathbb{N}$. Note that these are the only radii which make sense
  -- if $p^{-t} \leq r < p^{-t+1}$, then $B(x,r) = B(x,p^{-t})$.

  We take a moment to verify that the set  $\Bad(i,j)$ is indeed equal to
  the set $\Bad^*(\mathcal{R}, \beta, \rho_1, \rho_2)$. Fix $q \in \mathbb{Z}
  \setminus \{0\}$ and $(r_1,r_2) \in \mathbb{Z}^2$. Associated with
  the pair $((r_1, r_2),q)$ is the resonant point $R_{((r_1, r_2),q)}
  = (R_{(r_1,q)},R_{(r_2,q)})$. First, note that $|qx_s - r_s|_p =
  |q|_p \ d_1(x_s, R_{(r_s,q)})$ for $s \in \{1, 2\}$. However, $|q|_p
  \leq 1$ and so clearly $\Bad(i,j) \subseteq \Bad^*(\mathcal{R},
  \beta, \rho_1, \rho_2)$. Conversely, let $x \in \Bad^*(\mathcal{R},
  \beta, \rho_1, \rho_2)$.  We show that \eqref{psetup} is satisfied
  for $r$ and $q$.  If $(q,p) = 1$, then $|q|_p = 1$ and the
  inequality is immediate.  If $p^t | q$ for some $t \in \mathbb{N}$,
  but either $(r_1,p) = 1$ or $(r_2,p)=1$, the inequality is also
  satisfied. To see this, suppose that $(r_1,p) = 1$ and express
  $-r_1$ and $q x_1$ as power series in $p$. Clearly, the lowest
  exponent of $p$ in the expansion of $q x_1$ it at least $t$, whereas
  the expansion of $-r_1$ has a term with exponent zero. Hence the sum
  of the two must have a term of exponent zero, and so $|q x_1 -
  r_1|_p = 1$ and we are done. In the remaining case, when $p$ divides
  $q$, $r_1$ and $r_2$, we simply factor out the highest possible
  power of $p$ in the left hand side of \eqref{psetup} and the problem
  reduces to one of the previous cases. Thus, $\Bad^*(\mathcal{R},
  \beta, \rho_1, \rho_2) \subseteq \Bad(i,j)$.

  With reference to Theorem \ref{maingencorr}, the functions $\rho_1,
  \rho_2 $ satisfy condition (B*) and the measures $m_1:=\mu$ and
  $m_2:=\mu$ satisfy condition (A) with $\delta_1=\delta_2 = 1 $. We
  need to establish the existence of the collection $\cthfn$ where
  $F_n $ is an arbitrary closed rectangle of size $ 2k^{-n(1+i)}
  \times 2k^{-n(1+j)} $.  Here, we take $k=p^s$ and $\theta = p^{-t}$
  for some $s,t \in \mathbb{N}$ which will be chosen sufficiently
  large later on. In view of Lemma \ref{crucial}, there exists a
  disjoint collection $\cthfn$ of rectangles $2\theta F_{n+1} \subset
  \theta F_n $ such that \eqref{hgs1} of Theorem \ref{maingencorr} is
  satisfied. We now verify that \eqref{hgs2} is satisfied for any such
  collection.  This follows by modifying the `triangle' argument of
  \S\ref{basicgen} to the $p$-adic setting. So, let us assume that we
  have three resonant points (which by definition are rational points)
  $(r_1/q,r_2/q), (r'_1/q',r'_2/q') $ and $(r''_1/q'',r''_2/q'')$
  lying in some rectangle $\theta F_n$ with
  \begin{equation}
    k^n \leq \max_{1 \leq s,t,u \leq 2} \{ \vert r_s \vert, \vert r'_t
    \vert, \vert r''_u \vert , \vert q \vert \} < k^{n+1}.
    \label{pain}
  \end{equation}
  Suppose that they do not lie on a line. Then, they span a $p$-adic
  triangle $\Delta$. By results in Lutz \cite[Chapter I, \S
  4]{MR16:1003d}, the Haar measure $m$ of $\Delta$ is comparable to
  \begin{displaymath}
    \left\vert
      \det
      \begin{pmatrix}
        1 & r_1/q & r_2/q \\
        1 & r'_1/q' & r'_2/q' \\
        1 & r''_1/q'' & r_2''/q'' \\
      \end{pmatrix}
    \right\vert_p
    \neq 0.
  \end{displaymath}

  The determinant is a rational number with denominator $q q' q''$.
  As these are integers, the $p$-adic absolute value is $\leq 1$.
  Hence, the absolute value of the determinant is bounded below by the
  $p$-adic absolute value of the enumerator:
  \begin{equation*}
    N = r_1 r'_2 q'' - r_2 r'_1 q'' - r_1 q' r''_2 + r_2 q' r''_1 + q
    r'_1 r''_2 - q r'_2 r''_1.
  \end{equation*}
  This is an integer. In view of (\ref{pain}), we have that
  \begin{displaymath}
    |N| < 6 k^{3n+3}.
  \end{displaymath}
  We may assume without loss of generality that $N > 0$. Clearly, the
  $p$-adic valuation $v_p(N)$ (i.e. the number of times $p$ divides
  $N$) satisfies
  \begin{displaymath}
    v_p(N) < \log_p (6 k^{3n+3}).
  \end{displaymath}
  But $\vert N\vert_p = p^{-v_p(N)}$ so that
  \begin{displaymath}
    \vert N \vert_p > p^{-\log_p (6 k^{3n+3})} = 1/(6 k^{3n +3})  \ .
  \end{displaymath}
  Hence, there is a constant $C> 0$ such that $m(\Delta) > \ C/(6
  k^{3n+3}) $.  However, $\mu(\theta F_n) \leq \theta^2 k^{-3n}$ and
  on choosing $\theta^2 := p^{-2t} < C/(6 k^3)$ we obtain the desired
  contradiction; i.e. by choosing $t$ sufficiently large.  Thus, it
  there are two or more resonant points satisfying \eqref{pain} lying
  within $\theta F_n$ then they must lie on a $p$-adic line
  $\mathcal{L}$. It now follows that
  \begin{equation*}
    \text{l.h.s. of \eqref{hgs2}} \ \leq \ \#\{2\theta F_{n+1}
    \subset\cthfn: 2\theta F_{n+1} \cap \mathcal{L} \neq \emptyset\} \ .
  \end{equation*}
  A simple geometric argument, analogous to that employed in
  \S\ref{basicgen}, ensures that the line $\mathcal{L}$ can not pass
  through more than $C' \, k^{j+1}$ of the $2\theta F_{n+1}$
  rectangles. Here $C'>0$ is a constant independent of $k$.  On
  choosing $k:=p^s$ sufficiently large (i.e.  $s$ large enough), we
  ensure that $C' \, k^{1+j} < \kappa_1 \; k^{3} $ which establishes
  \eqref{hgs2} and thereby completes the proof of the theorem.
  \hfill $\Box$
\end{pf*}

Under suitable assumptions on subsets $\Omega_i$ of $\mathbb{Z}_p$
with measures satisfying condition (A), we can also obtain the
$p$-adic analogues of Theorems \ref{thm5} and \ref{thm6}. Of course,
to achieve this, one also needs to assume the natural $p$-adic
analogue of a measure being absolutely $\alpha$-decaying.

\subsection{Formal power series}

Apart from the $p$-adics, badly approximable elements have been
extensively studied over another locally compact ultra-metric field.
Let $\mathbb{F}$ be the finite field with $h$ elements. Thus, $h =
p^r$ for some prime $p$ and $r \in \mathbb{N} $. Now define
\begin{displaymath}
  \mathbb{F}((X^{-1})) := \left\{\sum_{i=-n}^\infty a_{-i} X^{-i} : n
    \in \mathbb{Z}, a_i \in \mathbb{F}, a_n \neq 0\right\} \cup
    \{0\}\ ,
\end{displaymath}
with an absolute value
\begin{displaymath}
  \left\Vert\sum_{i=-n}^\infty a_{-i} X^{-i}\right\Vert := h^n, \quad
  \left\Vert 0 \right\Vert := 0 \ .
\end{displaymath}
Under ordinary addition and multiplication, this is a locally compact
field. The closed unit ball $I = \{x \in \mathbb{F}((X^{-1})) : \Vert
x \Vert \leq 1\}$ is a compact subspace of this space.

In this section we consider the badly approximable analogue of
$\Bad(i_1,\ldots,i_N)$ in $I^N$.  Let $N \in \mathbb{N}$ and $i_1,
\ldots, i_N \geq 0$ such that $i_1 + \cdots + i_N = 1$.  Now define
the set $\Bad_{\mathbb{F}((X^{-1}))}(i_1,\dots, i_N)$ to consist of $x
:=(x_1, \ldots , x_N) \in \mathbb{F}((X^{-1}))^N$ for which there
exists a constant $c(x) > 0$ such that
\begin{displaymath}
  \max\{ \Vert q x_1 - p_1 \Vert^{1/i_1}, \dots, \Vert q x_N - p_N
  \Vert^{1/i_N} \} \geq c(x) \, \Vert q \Vert^{-1}
\end{displaymath}
for all $q, p_1, \dots, p_N \in \mathbb{F}\,[X]$ $(q \neq 0)$. Note
that in this setup, the polynomial ring $\mathbb{F}\, [X]$ plays the
role of the integers. When $i_1= \ldots = i_N= 1/N$, the corresponding
set will be denoted by $\Bad_{\mathbb{F}((X^{-1}))}(N)$. Niederreiter
and Vielhaber \cite{MR99f:94012} have shown that the set
$\Bad_{\mathbb{F}((X^{-1}))}(1)$ has full dimension. Using the
framework established in this paper, we are able to obtain the
complete result for the `weighted' simultaneous set.

\begin{thm}
  \begin{equation*}
    \dim \Bad_{\mathbb{F}((X^{-1}))} (i_1,\ldots, i_N) = N \ .
  \end{equation*}
  \label{thmformal}
\end{thm}

\begin{pf*}{Proof of Theorem \ref{thmformal} \emph{(Sketch)}. \ }
  As usual, we restrict our attention to the case $N = 2$ and write
  $\Bad_{\mathbb{F}((X^{-1}))}(i,j)$ for
  $\Bad_{\mathbb{F}((X^{-1}))}(i_1,i_2)$.  In view of the geometrical
  nature of our approach and the similarities between this situation
  and the preceding ones (in particular the $p$-adic case), we only
  outline the modifications needed to deal with the present situation
  in the briefest sense. The field $\mathbb{F}((X^{-1}))$ supports a
  Haar measure $m$ satisfying $m(B(c, h^{-t})) = h^{-t}$ for all $t
  \in \mathbb{Z}$. As was the case in the $p$-adics, these are the
  only balls for which a calculation is needed. Let $I$ denote the
  unit ball in this space.  We set $X_1 = X_2 = \mathbb{F}((X^{-1}))$,
  $\Omega_1 = \Omega_2 = I$ with the metrics induced by the absolute
  value and Haar measure defined above. We let $J = \{((p_1, p_2),q)
  \in \mathbb{F}[X]^2 \times \mathbb{F}[X] \setminus \{0\}\}$ and for
  any $((p_1, p_2),q) \in J$, we let $\beta_{((p_1, p_2),q)} = \Vert q
  \Vert$. The resonant sets $R_{((p_1, p_2),q)} = (p_1/q, p_2/q)$.
  Finally, define functions $\rho_1(x) = x^{-(i+1)}$ and $\rho_2(x) =
  x^{-(j+1)}$.  Clearly, the conditions of Theorem \ref{maingencorr}
  are satisfied and the set $\Bad_{\mathbb{F}((X^{-1}))}(i,j) \cap I^2
  = \Bad^*(\mathcal{R}, \beta, \rho_1,\rho_j)$.

  We establish the collection $\mathcal{C}(\theta F_n)$ by Lemma
  \ref{crucial}. The triangle argument works in this setting by
  results of Mahler \cite{MR2:350c} to calculate the measures of the
  sets involved. Note that in this case, the lower bound on the
  denominator is the important feature in the argument, so the proof
  differs from the $p$-adic case in this respect. Finally, maximal
  number of rectangles in $\mathcal{C}(2\theta F_{n+1})$ with
  non-trivial intersection with the resulting `line' is estimated by
  arguments as in the $p$-adic case. \hfill $\Box$
\end{pf*}

As in the $p$-adic setup, under appropriate assumptions we can also
obtain the formal power series analogues of Theorems \ref{thm5} and
\ref{thm6}. We have chosen to restrict ourselves to the simpler
situation, as this already yields new results and illustrates the
versatility of our framework.


\noindent\textbf{Acknowledgments: \, } SV would like to thank Dmitry
Kleinbock and Barak Weiss for the many e-mail conversations regarding
the development of our separate approaches. Also, he'd like to take
this opportunity to say `marra peeta' to Iona and Ayesha.

\end{document}